\documentclass{article}
\usepackage{xcolor}
\definecolor{PDHGblue}{RGB}{31,119,180}
\definecolor{Snaporange}{RGB}{230,120,20}
\usepackage[letterpaper,margin=1in]{geometry}
\usepackage{times}
\usepackage{natbib}

\usepackage[utf8]{inputenc}
\usepackage[T1]{fontenc}
\usepackage{amsmath,amssymb,amsthm}
\usepackage{booktabs}
\usepackage{graphicx}
\usepackage{algorithm}
\usepackage{algpseudocode}
\usepackage{longtable}
\usepackage{microtype}
\usepackage{pgfplots}
\usepackage{tikz}
\usetikzlibrary{arrows.meta}
\usepackage{xcolor}
\usepackage{hyperref}
\usepackage{url}

\theoremstyle{plain}
\newtheorem{theorem}{Theorem}
\newtheorem{lemma}{Lemma}
\newtheorem{proposition}{Proposition}
\newtheorem{corollary}{Corollary}
\theoremstyle{definition}
\newtheorem{assumption}{Assumption}
\theoremstyle{remark}
\newtheorem{remark}{Remark}

\newcommand{\proj}{\operatorname{proj}}
\newcommand{\E}{\mathbb{E}}
\newcommand{\Prob}{\mathbb{P}}
\newcommand{\Var}{\operatorname{Var}}
\newcommand{\CV}{\operatorname{CV}}
\newcommand{\Tdet}{T^{\mathrm{det}}}
\newcommand{\cT}{\mathcal{T}}
\newcommand{\cA}{\mathcal{A}}
\newcommand{\cS}{\mathcal{S}}

\title{First Passage Time in Geometric Resetting of\\ First-Order Methods in Linear Programming}

\author{%
  Yassine Koubaa$^{1,2,3}$ and Jakub Mare\v{c}ek$^{4}$\\[0.4em]
  \normalsize $^{1}$\textit{Artelys, Paris, France}\\
  \normalsize $^{2}$\textit{Universit\'e C\^ote d'Azur, Nice, France}\\
  \normalsize $^{3}$\textit{UiT The Arctic University of Norway, Troms\o, Norway}\\
  \normalsize $^{4}$\textit{Czech Technical University in Prague, Prague, Czech Republic}\\[0.3em]
    \normalsize \texttt{yassine.koubaa@artelys.com}\\   \normalsize \texttt{jakub.marecek@fel.cvut.cz}
}
\date{\normalsize (Dated: \today)}

\begin{document}

\maketitle

\begin{abstract}
Restarts are one of several ingredients of first-order methods for large linear programmes, such as the primal-dual hybrid gradient method (PDHG) that underlies PDLP, and the one with the most developed theory. We study the expected number of iterations such a method needs when it is restarted at random: at each iteration, with some probability, the iterate is replaced by a draw from a restart law. On an idealised family of instances with a plateau and a wall, randomised restarting beats PDLP by a factor that grows like the square root of the dimension. Once the method is placed in the framework of first passage under restart, the discrete-time renewal formula of Bonomo and Pal gives the expected iteration count under a constant restart probability, and its derivative at zero gives a dispersion criterion for infinitesimal rates. We derive a mean-residual-life characterisation of how geometric restarting affects the expected iteration count. Under independent restarts from a fixed distribution, we compare the expected number of iterations remaining after a run has not yet met the stopping criterion with the expected number needed from a fresh restart. If the former is never larger, restarting cannot reduce the expected iteration count; if it is never smaller, restarting cannot increase it. For base methods whose iteration count obeys a two-sided linear-rate bound, geometric restarting can improve the expected count only by a factor tending to one as the tolerance decreases.In a reduced-state control model with stochastically
monotone dynamics, the optimal state-dependent policy
is a two-level threshold rule with a non-degenerate
no-restart neighbourhood under an additional condition. On $68$ LP relaxations of MIPLIB 2017 models, a state-dependent stochastic rule in place of the adaptive restart of PDLP has a geometric-mean iteration ratio of $1.074$ on the $66$ instances on which both methods terminate, needs fewer iterations on $27$ of them, and needs one twentieth of the iterations on the instance on which PDLP is slowest.
\end{abstract}

\section{Introduction}

Linear programming on a large scale arises throughout operations research, machine learning and combinatorial optimisation. For three decades after Dantzig's formulation the simplex method had the field to itself \citep{dantzig1955generalized,fearnley2015complexity}; over the following four, interior-point methods rose to prominence \citep{karmarkar1984new,gill1986projected,andersen1996implementation,wright1997primal,gondzio2012interior}. In the past five years PDLP \citep{Applegate2021PracticalLL,applegate2026pdlp} and related first-order solvers \citep{lu2025cupdlp,lin2025pdcs,chen2026relationships} have become competitive on instances whose factorisation cannot be held in memory, or for which a moderate accuracy suffices. PDLP combines the primal-dual hybrid gradient (PDHG) iteration with presolve, diagonal preconditioning, adaptive step sizes, a primal weight and a restart scheme. Of these ingredients the restart has the most developed theory: \citet{Applegate_2022} prove that, under a sharpness condition, restarted PDHG converges linearly, whereas the guarantees for the unrestarted method are sublinear. Restarting too often discards the progress accumulated since the last restart; restarting too seldom leaves the iteration in a slow regime.

A question of the same shape has been studied in statistical physics under the name of stochastic resetting \citep{evans2011diffusion,Evans_2011,Evans_2014,Redner_2001,evans2020stochastic}. A diffusing particle seeks a target and is returned to its origin at random times. \citet{evans2011diffusion} showed that the return can shorten the mean time to reach the target, and \citet{reuveni2016optimal} and \citet{pal2017first} gave the criterion: an infinitesimal rate of resetting helps if and only if the coefficient of variation of the reset-free first-passage time exceeds one. An iterative optimisation method is, in outline, a process of the same kind: it has a target, the set of points at which the termination test is satisfied, and it may be sent to a fresh point at random times. Between restarts, however, its motion is deterministic, and its only randomness is the restart itself. This paper asks when randomised restarting reduces the expected iteration count of a first-order method for linear programming, and how the restart probability should depend on the state of the iteration. Section~\ref{sec:wall} gives an idealised family on which randomised restarting wins by a margin that grows with the dimension; it motivates what follows. Our contributions are the following.
\begin{enumerate}
\item \emph{A renewal model of restarted first-order methods.} We cast a first-order method with randomised restarts as a first passage under restart. The discrete-time renewal formula and dispersion criterion of \citet{bonomo2021first} then give the expected iteration count under a constant restart probability (Theorem~\ref{thm:renewal}, Corollary~\ref{cor:cv}).
\item \emph{A global analysis. }We represent $T_p-T_0$ through the mean residual life of the restart-free count, where $T_p$ is the expected iteration count under geometric restarting with probability $p$ and $T_0$ that without restarting (Theorem~\ref{thm:mrl}). If an unfinished run never expects to need more (fewer) further iterations than a fresh one, restarting cannot help (hurt) (Corollary~\ref{cor:nbue}). Under a two-sided linear-rate bound the best speed-up tends to one as the tolerance decreases, and sufficiently small restart probabilities hurt (Proposition~\ref{prop:bounded}, Corollary~\ref{cor:linear}).
\item \emph{An optimal state-dependent policy.} In a reduced-state control model, we prove that a two-level
threshold restart policy is optimal under a stochastic monotonicity assumption (Theorem~\ref{thm:threshold}).
This assumption is verified for a diffusion surrogate,
but remains unverified for projected PDHG.
Under an additional condition, the optimal policy
never restarts near the target (Corollary~\ref{cor:plow}).

\item \emph{Experiments.}
A heuristic state-dependent stochastic restart rule
performs comparably to PDLP on $66$ MIPLIB relaxations,
with one substantial instance-level gain.
On \texttt{danoint}, the renewal formula
correctly predicts the direction of the observed
slowdown under constant-rate restarting
(Section~\ref{sec:experiments}).
On coupled-wall PageRank LPs, sPDLP with snapping
requires up to $13.5$ times fewer iterations than PDLP
with fixed steps and no diagonal preconditioning
(Figure~\ref{fig:pagerank-iterations},
Appendix~\ref{app:pagerank}).
Across $18$ instances, it also achieves the lowest
penalised shifted geometric-mean solve time among
the evaluated solvers under the stated penalty convention.
\end{enumerate}

\subsection{A motivating example: a plateau and a wall}
\label{sec:wall}

The following idealised family of LP instances shows that randomised restarting can beat PDHG with PDLP-style averaging restarts by a factor that grows with the dimension, and shows what such a gain requires. Take a variable $w\in[0,1]$ with cost $-\delta$, together with an independent LP subproblem (called the block), $\min\{\hat c^\top\hat x:\hat K\hat x\ge\hat h,\ \hat l\le\hat x\le\hat u\}$, whose constraint matrix has $\|\hat K\|_2=\sqrt n$. Since $w$ does not appear in the block's constraints, its optimum is $w=1$. Run PDHG with $\tau=\sigma$ and $\tau\|\hat K\|_2<1$, as the step-size condition requires. One step moves $w$ to $\min\{1,w+\tau\delta\}$: $w$ crosses a plateau at speed $\tau\delta$ and stops at the wall $w=1$. Write $z^k=(w^k,\hat x^k,\hat y^k)$ for the iterate, the primal-dual point after $k$ iterations. The target is the set of points with $w=1$ and $(\hat x,\hat y)\in\mathcal B$, where $\mathcal B$ is a prescribed target set for the block, such as that of a stopping criterion.

\begin{theorem}[Plateau and wall]\label{thm:wall}
Suppose that from every starting point the block's primal-dual iterates reach its target set $\mathcal B$ within $f\ge2$ PDHG iterations independently of $w$, and let $L:=\lceil 1/(\tau\delta)\rceil$, so that $L\ge\sqrt n/\delta$.
\begin{enumerate}
\item[(i)] Any method whose iterates are PDHG steps or convex combinations of earlier iterates (weighted averages with
nonnegative weights summing to one) needs at least $L$ iterations from $w=0$. This includes PDHG with PDLP-style averaging restarts
and the unsnapped stochastic restart rule of
Section~\ref{sec:policy}, whose restart destinations
lie between the current iterate and its running
average.
\item[(ii)] Let a restart draw $w\in\{0,1\}$ with equal probabilities and the block from any fixed law. Geometric restarting at rate $p=1/f$ reaches the target in at most $2\mathrm{e}f$ iterations in expectation.
\end{enumerate}
\end{theorem}

The proof is in Appendix~\ref{app:proof_wall}. The gain comes from a restart law that places mass on the wall itself, so that a restarted run skips the plateau with probability one half. Restarts confined to the convex hull of the trajectory cannot skip it, which is the content of part~(i).

\begin{figure}[H]
\centering
\begin{tikzpicture}
\begin{axis}[
    width=0.91\linewidth,
    height=6.1cm,
    xmin=0.5,xmax=10.6,
    ymin=0,ymax=470,
    xtick={1,3,10},
    xticklabels={$10^4$,$3\cdot 10^4$,$10^5$},
    xlabel={Number of nodes $N$},
    ylabel={Iterations ($10^3$)},
    grid=major,
    grid style={gray!20},
    legend style={font=\footnotesize,at={(0.02,0.98)},anchor=north west,
                  draw=none,fill=white,fill opacity=0.85,
                  text opacity=1},
    tick label style={font=\footnotesize},
    label style={font=\small},
    every axis plot/.append style={line width=0.9pt}
]
\addplot[PDHGblue,mark=square*,solid,mark size=2pt] coordinates
  {(1,118.120) (3,223.680) (10,385.000)};
\addlegendentry{PDLP}

\addplot[PDHGblue,mark=*,densely dotted,mark size=2pt] coordinates
  {(1,150.560) (3,236.280) (10,415.400)};
\addlegendentry{sPDLP}

\addplot[Snaporange,mark=triangle*,dashed,mark size=2.4pt,line width=1.4pt] coordinates
  {(1,14.560) (3,21.742) (10,28.591)};
\addlegendentry{sPDLP with snapping}

\addplot[gray!70,mark=none,densely dotted,domain=1:10,samples=80]
  {125*sqrt(x)};
\addlegendentry{$1.25\sqrt{N}/\delta$}
\end{axis}
\end{tikzpicture}
\caption{Iterations on two-community PageRank LPs with one coupled wall ($k=1$), at the three sizes on which all three PDLP variants terminate; fixed step sizes, no diagonal preconditioning. The dotted grey curve is the plateau length $1/(\tau\delta)\approx1.25\sqrt N/\delta$ suggested by Theorem~\ref{thm:wall}. With snapping every run ends at $w_1=1$; without it, a run may meet the tolerance with $w_1<1$. Data from Table~\ref{tab:c4-coupled-walls}.}
\label{fig:pagerank-iterations}
\end{figure}
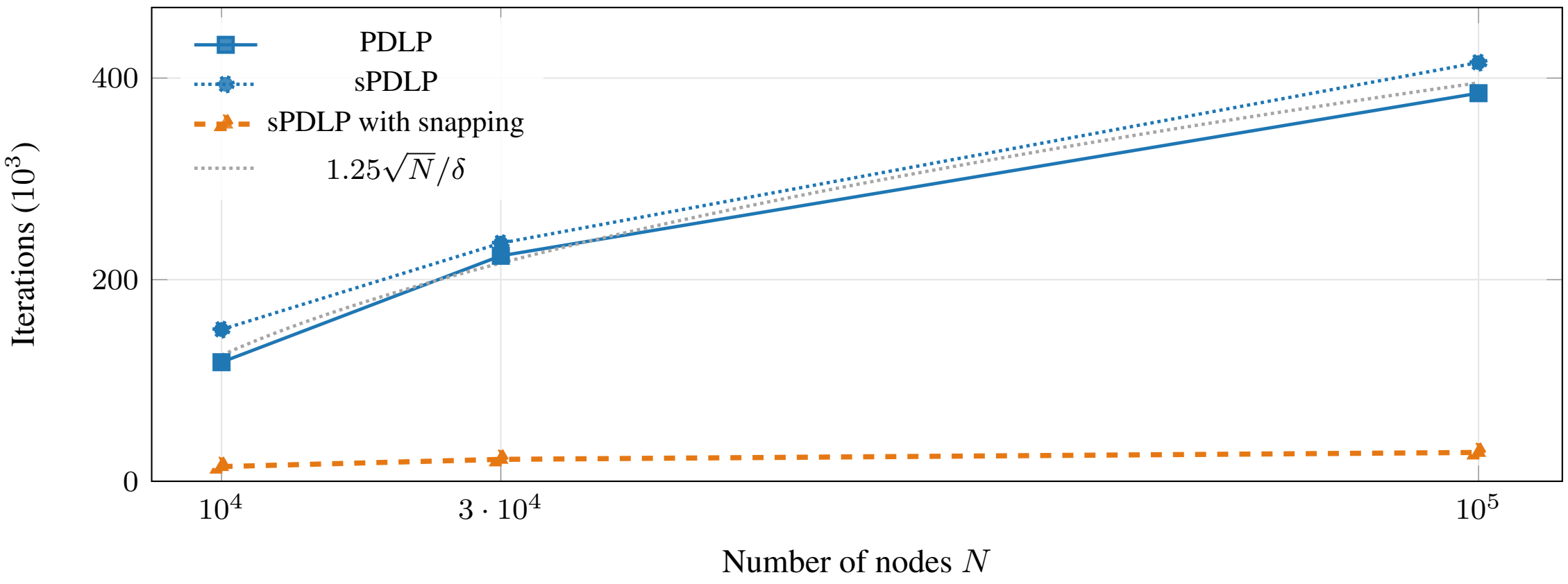

Appendix~\ref{app:plateau-experiments} gives a numerical illustration and the limitations of the construction. Snapping lets a restart move selected coordinates to their bounds instead of averaging earlier iterates (Section~\ref{sec:policy}). A wall variable that appears in no constraint would, however, be removed by presolve. To test whether snapping remains
effective when walls enter the constraints,
we therefore consider two-community PageRank LPs,
in which each wall variable appears in two
PageRank equations. Figure~\ref{fig:pagerank-iterations} shows the outcome; Appendix~\ref{app:pagerank} gives the construction and the full results.

\section{Related work}

\citet{Applegate_2022} gave the first rigorous analysis of restarted PDHG for linear programming, proving linear convergence under sharpness for several restart schemes; \citet{applegate2026pdlp} document the engineering, and \citet{applegate2024infeasibility} extend the analysis to infeasibility detection. Later work refines the rate through the geometry of the problem \citep{lu2024geometry,xiong2024role}, gives guarantees in terms of limiting error ratios and of the uniqueness of the optimiser \citep{xiong2026computational,xiong2026accessible}, accelerates by restarting a Halpern iteration \citep{lu2024restarted}, and carries the method to quadratic programming \citep{lu2025practical} and to GPUs \citep{lu2025cupdlp,lin2025pdcs,chen2026relationships}. For accelerated gradient methods, adaptive restarts are due to \citet{odonoghue2015adaptive}; \citet{roulet2020sharpness}, \citet{fercoq2019adaptive}, \citet{renegar2022simple} and \citet{adcock2025restarts} give schedules that are near-optimal under sharpness or quadratic growth. In all of these the restart is a deterministic function of the trajectory. That randomised restarts shorten heavy-tailed running times is classical in combinatorial search: \citet{luby1993optimal} give optimal universal schedules for Las Vegas algorithms, and \citet{gomes1998boosting} and \citet{gomes2000heavy} documented heavy-tailed running times of backtrack search for satisfiability and constraint satisfaction and showed that randomised restarts remove them. The bimodal example of Remark~\ref{rem:notnecessary} is of this type.

\citet{evans2011diffusion,Evans_2011,Evans_2014} introduced diffusion with stochastic resetting; \citet{evans2020stochastic} survey the field. \citet{reuveni2016optimal} and \citet{pal2017first} gave the criterion for infinitesimal rates, \citet{pal2019landau} analyse the transitions of the optimal rate, \citet{pal2022inspection} read the criterion through the inspection paradox, and \citet{eliazar2020mean} analyse sharp restart in continuous time through the residual lifetime and the hazard rate of the completion time; see also \citet{Ahmad_2019,aspman2026first}. \citet{bonomo2021first} set up first passage under restart in discrete time for a generic process whose memory is erased at each restart, derive the mean completion time under geometric restart and the criterion $\CV^2>1-1/\E[N]$, and apply them to random walks on a one-dimensional lattice. Theorem~\ref{thm:renewal} and Corollary~\ref{cor:cv} below are their eq.~(20) and criterion~(25), up to the convention for the timing of a restart; we give self-contained proofs because the later results are built from the same decomposition.
While often seen as restricted to univariate settings, the renewal framework of \citet{pal2017first} and \citet{bonomo2021first} is not confined to one dimension, since it depends only on the law of the scalar reset-free time. What the resetting literature has in common is its object: a stochastic process, typically Brownian motion or a lattice walk in one dimension or in a few \citep{Evans_2014}, returned to a fixed point or to a prescribed law. Our setting differs in three respects. The motion between restarts is a deterministic iteration in dimension $n$ plus the number of constraints, so that all randomness of the restart-free count is induced by the restart destination. The restart law is chosen by the algorithm and depends on the trajectory, which is why the renewal assumption is an idealisation here. And the quantities that matter are tied to the convergence rate of the base method. To our knowledge, Theorem~\ref{thm:wall}, the representation of Theorem~\ref{thm:mrl} at every rate and its dichotomy for geometric restart in discrete time, Proposition~\ref{prop:bounded} and Corollary~\ref{cor:linear}, and the treatment of state-dependent restarting of an optimisation method as a monotone dynamic programme in Section~\ref{sec:threshold}, whose argument is itself standard \citep{puterman1994markov}, are new.

\section{Setting}
\label{sec:setting}

We consider the linear programme
\begin{equation}
  \min_{x\in\mathbb{R}^n}\ c^\top x \quad \text{s.t.}\ Gx\ge h,\ Ax=b,\ l\le x\le u,
\label{eq:LP}
\end{equation}
which, with $K^\top=(G^\top,A^\top)$ and $q^\top=(h^\top,b^\top)$, is the saddle-point problem $\min_{x\in X}\max_{y\in Y}\ c^\top x - y^\top Kx + q^\top y$ over $X=\{x:l\le x\le u\}$ and $Y=\{y=(y_G,y_A):y_G\ge 0\}$. PDLP applies the alternating projected steps
\begin{equation}
x^{k+1} = \proj_X\!\bigl(x^k - \tau(c - K^\top y^k)\bigr), \qquad
y^{k+1} = \proj_Y\!\bigl(y^k + \sigma(q - K(2x^{k+1}-x^k))\bigr),
\label{eq:pdlp}
\end{equation}
with step sizes $\tau,\sigma>0$, adapted during the run, to the iterates $z^k=(x^k,y^k)$, and restarts to their running average $\bar z^k$. Progress is monitored, as in \citet{Applegate2021PracticalLL}, through the convergence score
\begin{equation}
s^k := r_{\mathrm{p}}^k + r_{\mathrm{d}}^k + g^k,
\label{eq:score}
\end{equation}
the sum of the relative primal residual, the relative dual residual and the relative gap, defined in Appendix~\ref{app:score}. We write $s(z)$ for the score of an arbitrary point $z$; the solver stops when $\bar s^k:=s(\bar z^k)\le\varepsilon_{\mathrm{opt}}$.

Let $\Psi$ denote one iteration of a deterministic base method. We model the restarted trajectory as
\begin{equation}
z^{k+1} =
\begin{cases}
  \Psi(z^k), & \text{with probability } 1 - p_k, \\
  \tilde{z}^{k+1}\sim\mu_k(\cdot\mid z^k), & \text{with probability } p_k,
\end{cases}
\label{eq:process}
\end{equation}
where $p_k\in[0,1]$ is the restart probability and $\mu_k$ a restart law. The analysis applies to any deterministic base method. In the rule of Section~\ref{sec:policy} the stochastic restart replaces the adaptive restart of PDLP, so there $\Psi$ is the PDLP iteration with its own restart switched off. To obtain a tractable one-dimensional state we project $z^k$ onto the \emph{progress variable}

\begin{equation}
R_k :=
\frac{1}{\tau}\|x^{k+1}-x^k\|_2^2
+\frac{1}{\sigma}\|y^{k+1}-y^k\|_2^2
\end{equation}
a by-product of every step, which tends to zero as the iterates converge. For the empirical diagnostic in
Appendix~\ref{app:table}, we evaluate this quantity using
the iteration-dependent step sizes
$\tau_k=\eta_k/\omega_k$ and
$\sigma_k=\eta_k\omega_k$ in the solver's scaled space. For the reduced model, the target is $\cT_\varepsilon=\{R\le\varepsilon\}$, and the mean first passage time from $R$ is $T(R)=\E[\inf\{k:R_k\le\varepsilon\}\mid R_0=R]$. A diffusion surrogate for $R_k$, used only to exhibit a model in which the monotonicity hypothesis of Section~\ref{sec:threshold} holds, is described in Appendix~\ref{app:surrogate}.

\section{Geometric restarting: an exact analysis}
\label{sec:exact}

\subsection{Standing assumptions}
\label{sec:exact_setting}

In this section the restart probability is a constant, $p_k\equiv p$, and the restart law is fixed. Fix a target set $\cT$; the reader may take $\cT_\varepsilon$ or the termination test $\{s\le\varepsilon_{\mathrm{opt}}\}$, as nothing below depends on the choice. Write $\mu$ for the law of the point produced by a restart.

\begin{assumption}[Restart-free iteration count]\label{ass:N}
Let $N$ be the number of iterations the base method needs to reach $\cT$ when started from a point drawn from $\mu$ and left to run without further restarts. We assume that $N$ takes values in $\{1,2,\dots\}\cup\{\infty\}$, that $\Prob[N<\infty]>0$, and that successive restart epochs produce independent copies $N_1,N_2,\dots$ of $N$, independent also of the restart coins.
\end{assumption}

Since the base map is deterministic, all the randomness in $N$ comes from the restart: $N=N(\tilde z)$ with $\tilde z\sim\mu$. The value $N=\infty$ covers stalling. Independence holds when each restart draws afresh from a fixed law, as in Theorem~\ref{thm:wall}; for the trajectory-dependent restart action of Section~\ref{sec:policy} it fails, and the results describe a renewal idealisation whose error we do not quantify. A restart is not charged as an iteration; Remark~\ref{rem:convention} in the appendix gives the constants under the opposite convention.

\subsection{The renewal formula}
\label{sec:renewal}

\begin{theorem}[Mean first passage time under geometric restarting; \citealp{bonomo2021first}]\label{thm:renewal}
Let Assumption~\ref{ass:N} hold, let $p\in(0,1)$, and let $T_p$ be the expected number of iterations to reach $\cT$ when the restart is applied, independently at each iteration, with probability $p$, starting from a fresh draw from $\mu$. Write
\[
  \phi(p) \;:=\; \E\bigl[(1-p)^{N}\bigr]
  \;=\;\sum_{t\ge 1}\Prob[N=t]\,(1-p)^{t}
  \qquad(\text{with }(1-p)^{\infty}:=0).
\]
Then $\phi(p)\in(0,1)$ and
\begin{equation}
  T_p \;=\; \frac{(1-p)\bigl(1-\phi(p)\bigr)}{p\,\phi(p)}\;<\;\infty .
  \label{eq:renewal}
\end{equation}
Moreover $T_p\to\E[N]$ as $p\downarrow 0$, where $\E[N]=\infty$ is permitted on both sides.
\end{theorem}

The proof (Appendix~\ref{app:proof_renewal}) decomposes the run into epochs. Let $\Gamma$ be the geometric index of the first restart. An epoch ends in absorption if and only if $\Gamma\ge N$, an event of probability $\E[(1-p)^{N-1}]=\phi(p)/(1-p)$, and it consumes $\min(N,\Gamma)$ iterations in either case; Wald's identity then gives $T_p=\E[\min(N,\Gamma)]/\Prob[\Gamma\ge N]$.

\begin{corollary}[Restarting removes stalling]\label{cor:finite}
If $\Prob[N=\infty]>0$ then $T_0=\E[N]=\infty$, whereas $T_p<\infty$ for every $p\in(0,1)$, since $\phi(p)>0$.
\end{corollary}

This is the discrete analogue of the observation of \citet{Evans_2011} that resetting renders an infinite mean first passage time finite. It asserts only finiteness: by Proposition~\ref{prop:bounded}, $T_p$ grows exponentially in the number of iterations a concentrated $N$ requires.

\begin{corollary}[Dispersion criterion; \citealp{bonomo2021first}]\label{cor:cv}
Suppose in addition that $\Prob[N<\infty]=1$ and $\E[N^2]<\infty$. Then $p\mapsto T_p$ is right-differentiable at $0$, with
\begin{equation}
  \frac{d T_p}{dp}\Big|_{p=0^{+}}
  \;=\;\tfrac12\Bigl(\E[N]^2-\E[N]-\Var(N)\Bigr),
  \label{eq:cvderiv}
\end{equation}
so that, writing $\CV(N):=\sqrt{\Var(N)}/\E[N]$,
\begin{equation}
  \exists\,p>0:\ T_p<T_0
  \qquad\Longleftarrow\qquad
  \CV(N)^2>1-\frac{1}{\E[N]} .
  \label{eq:cvcriterion}
\end{equation}
If the inequality in~\eqref{eq:cvcriterion} is reversed and strict, then $T_p>T_0$ for all sufficiently small $p>0$.
\end{corollary}

The proof is in Appendix~\ref{app:proof_dispersion}. As $\E[N]\to\infty$ the criterion becomes the $\CV>1$ criterion of \citet{reuveni2016optimal} and \citet{pal2017first}; the correction $-1/\E[N]$ reflects the granularity of the iteration. A deterministic count $N\equiv t$ gives $T_p=\Tdet_p(t):=\sum_{k=0}^{t-1}(1-p)^{-k}$, increasing in $p$, and a geometric count gives $T_p=\E[N]$ for every $p$ (Remark~\ref{rem:extremes}).

\begin{remark}[The criterion is sufficient only]\label{rem:notnecessary}
Fix an integer $t>1$ and let $N=1$ or $N=t$, each with probability $\tfrac12$. Then $\E[N]=\tfrac{1+t}2$ and $\Var(N)=\tfrac{(t-1)^2}4\le\E[N]^2-\E[N]$, so criterion~\eqref{eq:cvcriterion} fails and small restart probabilities hurt. Yet letting $p\uparrow1$ in~\eqref{eq:renewal}, with $\phi(p)\sim\tfrac12(1-p)$, gives $T_p\to2$, whereas $T_0=\tfrac{1+t}2$. For $t=101$, the derivative at zero is $+25$ while $\inf_pT_p=2$. This is the mechanism of Las Vegas restart schedules \citep{luby1993optimal}, of Theorem~\ref{thm:wall}, and in the resetting literature of the discontinuous transition of \citet{pal2019landau}: when the running time has a fast mode, one should re-draw aggressively until it is hit, and a derivative at $p=0$ cannot see this.
\end{remark}

\subsection{The sign of the gain at every rate}
\label{sec:mrl}

For a given $p$, is $T_p$ smaller or larger than $T_0$? The answer is governed by the \emph{mean residual life} of $N$, $m(j):=\E[N-j\mid N>j]$ for integers $j\ge0$ with $\Prob[N>j]>0$, the expected number of further iterations needed once $j$ have been spent without success; write $m:=m(0)=\E[N]$. Restarting after $j$ iterations exchanges $m(j)$ for $m$. The following identity says that $T_p-T_0$ is a geometrically weighted average of such exchanges.

\begin{theorem}[Mean-residual-life representation]\label{thm:mrl}
Let Assumption~\ref{ass:N} hold with $\Prob[N<\infty]=1$ and $m<\infty$. Then for every $p\in(0,1)$,
\begin{equation}
  T_p-T_0
  \;=\;\frac{1-p}{\phi(p)}\sum_{j\ge1}p(1-p)^{j-1}\,\Prob[N>j]\,\bigl(m-m(j)\bigr).
  \label{eq:mrl}
\end{equation}
\end{theorem}

The proof (Appendix~\ref{app:proof_mrl}) takes three lines from $T_p=\E[\min(N,\Gamma)]/\Prob[\Gamma\ge N]$ and $\E[N]=\E[\min(N,\Gamma)]+\E[(N-\Gamma)^+]$. Every weight in~\eqref{eq:mrl} is positive, so the sign of $T_p-T_0$ is decided by the signs of $m-m(j)$.

\begin{corollary}[Global dichotomy]\label{cor:nbue}
Under the hypotheses of Theorem~\ref{thm:mrl}:
\begin{enumerate}
  \item[(i)] if $m(j)\le m$ for every $j\ge1$ with $\Prob[N>j]>0$, that is, if $N$ is \emph{new better than used in expectation} (NBUE), then $T_p\ge T_0$ for every $p\in(0,1)$, with strict inequality at every $p$ unless $N$ is geometric;
  \item[(ii)] if $m(j)\ge m$ for every such $j$ (NWUE), then $T_p\le T_0$ for every $p\in(0,1)$;
  \item[(iii)] if the discrete hazard $H(t):=\Prob[N=t\mid N\ge t]$, with the convention $H(t)=1$ whenever
$\mathbb P[N\ge t]=0$, is non-decreasing in $t$ (increasing failure rate, IFR), then $N$ is NBUE and (i) applies.
\end{enumerate}
\end{corollary}

The proof is in Appendix~\ref{app:proof_mrl}. Statement (iii) is the algorithmic reading: if a restart-free run that has gone on longer without terminating is thereby more likely to terminate at the next iteration, then no constant restart rate can help. It is the opposite of the heavy-tailed situation of Remark~\ref{rem:notnecessary}, in which a long run signals a slow mode that should be abandoned. Expanding~\eqref{eq:mrl} at $p=0$ recovers~\eqref{eq:cvderiv}. When the restart-free count is concentrated, bounds are available at every rate.

\begin{proposition}[Bounded dispersion]\label{prop:bounded}
Let Assumption~\ref{ass:N} hold with $t_-\le N\le t_+$ almost surely, $1\le t_-\le t_+<\infty$. Then for every $p\in(0,1)$, with $\Tdet_p$ as above,
\begin{equation}
  \Tdet_p(t_-)\le T_p\le \Tdet_p(t_+),
  \quad\text{hence}\quad
  T_p\ge t_-\ \text{ and }\
  \frac{T_0}{\inf_{p}T_p}\le\frac{\E[N]}{t_-}\le\frac{t_+}{t_-}.
  \label{eq:bounded}
\end{equation}
If moreover $(t_+-t_-)^2<4\,t_-(t_--1)$, then criterion~\eqref{eq:cvcriterion} fails strictly and $T_p>T_0$ for all sufficiently small $p>0$.
\end{proposition}

The proof (Appendix~\ref{app:proof_bounded}) sandwiches $\min(N,\Gamma)$ and $\Prob[\Gamma\ge N]$ between their values at $t_\pm$ and uses Popoviciu's bound on the variance. The condition holds when the slowest restart point needs fewer than about three times the iterations of the fastest. Note also that $\Tdet_p(t_-)\ge(1-p)^{-(t_--1)}$: at a fixed rate the expected count grows exponentially in $t_-$.

\begin{corollary}[Linearly convergent base methods]\label{cor:linear}
Suppose there are constants $\gamma\in(0,1)$, $0<\rho_-\le\rho_+<\infty$ and $t_0\in\mathbb{Z}_{\ge 0}$ such that, for $\mu$-almost every restart point, the restart-free count $N(\varepsilon)$ to tolerance $\varepsilon$ satisfies
\begin{equation}
  t_-(\varepsilon):=\Bigl\lceil\frac{\log(\rho_-/\varepsilon)}{\log(1/\gamma)}\Bigr\rceil
  \;\le\;N(\varepsilon)\;\le\;
  \Bigl\lceil\frac{\log(\rho_+/\varepsilon)}{\log(1/\gamma)}\Bigr\rceil+t_0=:t_+(\varepsilon)
  \qquad\text{for all }\varepsilon\in(0,\rho_-).
  \label{eq:linear_bracket}
\end{equation}
Then, as $\varepsilon\downarrow0$: (i) the best speed-up available from geometric restarting at any rate satisfies $T_0(\varepsilon)/\inf_pT_p(\varepsilon)\le t_+(\varepsilon)/t_-(\varepsilon)\to1$; (ii) there is $\varepsilon_0>0$ such that for every $\varepsilon<\varepsilon_0$, $T_p(\varepsilon)>T_0(\varepsilon)$ for all sufficiently small $p>0$.
\end{corollary}

The proof is in Appendix~\ref{app:proof_bounded}. An upper bound of the form in~\eqref{eq:linear_bracket} is what the sharpness analysis of \citet{Applegate_2022} provides for restarted PDHG, up to constants, with $\rho$ the distance of the restart point from the optimal set; the lower bound is an assumption that the method converges no faster than linearly. The spread of log-distances of restart points then becomes a spread of iteration counts that is $O(1)$ while the counts themselves grow like $\log(1/\varepsilon)$. In the family of Theorem~\ref{thm:wall} the plateau makes $t_0$ grow with $n$, which is why the corollary does not apply there.

\subsection{Using the formula}
\label{sec:practice}

Theorem~\ref{thm:renewal} turns the choice of a constant rate into a measurement of the whole curve. Draw $B$ restart points $\tilde z^{(1)},\dots,\tilde z^{(B)}\sim\mu$, run the restart-free method from each under a fixed budget, record $N_1,\dots,N_B$, and form the plug-in curve
\begin{equation}
  \hat\phi(p)=\frac1B\sum_{i=1}^{B}(1-p)^{N_i},
  \qquad
  \hat T_p=\frac{(1-p)\bigl(1-\hat\phi(p)\bigr)}{p\,\hat\phi(p)}.
  \label{eq:plugin}
\end{equation}
Its minimiser over a grid of $p$ estimates the best constant rate. Runs that exhaust the budget contribute $N_i=\infty$, yielding a conservative upper estimate of the uncensored plug-in value, which may substantially overestimate it for small $p$ and distort the estimated optimal rate. When all pilot runs terminate, the same data give an empirical mean residual life, whose profile provides a diagnostic for the conditions of Corollary~\ref{cor:nbue}.

\section{State-dependent restarting: the threshold policy}
\label{sec:threshold}

We now let the restart probability depend on the state. Let $\cS\subseteq[0,\infty)$ be a Borel state space carrying the reduced state $R$, with the usual order; let $\cT=\{R\le\varepsilon\}$; let $P(\,\cdot\mid R)$ be the transition kernel of one iteration of the base method in the reduced state, which presupposes a Markov approximation of $R_k$; and let $\mu$ be the law of the reduced state after a restart. A \emph{policy} is a measurable map $\pi:\cS\to\cA$, where $\cA\subseteq[0,1]$ is a non-empty compact set of admissible restart probabilities, and the process evolves as~\eqref{eq:process} with $p_k=\pi(R_k)$. Write $T_\pi(R)$ for the expected number of iterations to reach $\cT$ from $R$ under $\pi$, and $T^\star:=\inf_\pi T_\pi$. This is a dynamic programme with non-negative costs in which a restart is an action with transition kernel $\mu$; its Bellman recursion lets a restart pre-empt an absorption, a convention that differs from that of Section~\ref{sec:renewal} by a factor $1-p$ and affects no statement below (Remark~\ref{rem:preempt}).

\begin{assumption}[Regularity]\label{ass:ssp}
$\cA$ is compact with $\max\cA>0$; the kernels $P$ and $\mu$ are Borel; and some policy $\pi$ has $T_\pi$ bounded on $\cS$.
\end{assumption}

Under Assumption~\ref{ass:ssp}, let $\bar T>0$ be a
finite uniform upper bound on $T^\star$.
If Assumption~\ref{ass:N} holds for the reduced process
and $\cA\cap(0,1)\ne\emptyset$, an explicit choice of
such a bound is provided by the constant policy
$\pi\equiv\bar a$, for any $\bar a\in\cA\cap(0,1)$,
which has, by Remark~\ref{rem:preempt},
\begin{equation}
  T_\pi\;\le\;\bar T\;=\;\frac1{\bar a}
  +\frac{T_{\bar a}}{1-\bar a}\;<\;\infty
  \label{eq:Tbar}
\end{equation}
uniformly on $\cS$, with $T_{\bar a}$ finite by
Theorem~\ref{thm:renewal}; hence $T^\star\le\bar T$.
Since every one-stage cost equals $1$, value iteration
from the zero function increases to $T^\star$,
$T^\star$ satisfies the Bellman equation, and a
stationary policy attaining the minimum in it is
optimal \citep[Ch.~8--9]{bertsekas1978stochastic}.

\begin{lemma}[Bellman equation and bang--bang structure]\label{lem:bangbang}
Under Assumption~\ref{ass:ssp}, $T^\star\equiv0$ on $\cT$ and, for $R\notin\cT$,
\begin{equation}
  T^\star(R)
  \;=\;1+\min_{a\in\cA}
  \Bigl[(1-a)\,Q^\star(R)+a\,C^\star\Bigr],
  \qquad
  Q^\star(R):=\!\int\! T^\star\,dP(\cdot\mid R),
  \quad
  C^\star:=\!\int\! T^\star\,d\mu .
  \label{eq:bellman}
\end{equation}
The bracket is affine in $a$, so the action $\max\cA$ is optimal where $Q^\star(R)>C^\star$, the action $\min\cA$ where $Q^\star(R)<C^\star$, and every action where they are equal.
\end{lemma}

One should therefore restart aggressively exactly when the expected remaining time from where the iteration is heading exceeds that from a fresh restart. To conclude that the switching set is an interval we need the reduced dynamics to be monotone.

\begin{assumption}[Stochastic monotonicity]\label{ass:mono}
If $R\le R'$ then $P(\cdot\mid R)$ is stochastically dominated by $P(\cdot\mid R')$: $\int f\,dP(\cdot\mid R)\le\int f\,dP(\cdot\mid R')$ for every bounded non-decreasing $f$.
\end{assumption}

Lemma~\ref{lem:diffmono} in Appendix~\ref{app:surrogate} verifies this for the Euler--Maruyama discretisation of the diffusion surrogate. It is not verified for the kernel induced on $R_k$ by one projected PDHG step~\eqref{eq:pdlp}, whose projections onto a box and a cone are non-expansive but not obviously monotone once composed with the norm defining $R_k$; we know neither a proof nor a counterexample. The theorem below is thus exact for the surrogate, and only suggestive for the algorithm. An empirical diagnostic on \texttt{danoint}
(Appendix~\ref{app:danoint-monotonicity}) shows small violations of stochastic
monotonicity over one iteration, but substantially
larger violations over 40 iterations, the interval
between restart checks. These observations neither establish nor
disprove Assumption~\ref{ass:mono} for projected PDHG.

\begin{theorem}[The two-level threshold policy is optimal]\label{thm:threshold}
Let Assumptions~\ref{ass:ssp} and~\ref{ass:mono} hold and let $\cT=\{R\le\varepsilon\}$. Then:
\begin{enumerate}
  \item[(i)] $T^\star$ is non-decreasing on $\cS$, and so is $Q^\star$;
  \item[(ii)] the set $\{R:Q^\star(R)>C^\star\}$ is an up-set, of the form 
    $\cS\cap(\beta^\star,\infty)$ or
    $\cS\cap[\beta^\star,\infty)$, where
    \begin{equation}
      \beta^\star:=\inf\{R\in\cS: Q^\star(R)>C^\star\}\in[0,\infty];
      \label{eq:betastar}
    \end{equation}
  \item[(iii)] consequently the stationary policy
    \begin{equation}
     \pi^\star(R)=
     \begin{cases}
       \max\cA, & Q^\star(R)>C^\star,\\
       \min\cA, & Q^\star(R)\le C^\star,
     \end{cases}
      \label{eq:optpolicy}
    \end{equation}
    is optimal.
\end{enumerate}
In particular, $\cA=[0,\bar p]$ gives an optimal policy that does not restart for $R<\beta^\star$.
\end{theorem}
\par\smallskip
\noindent
The proofs of the lemma and the theorem are in Appendix~\ref{app:proof_threshold}.

\begin{corollary}[No restarting near convergence]\label{cor:plow}
If $0\in\cA$, then the optimal policy~\eqref{eq:optpolicy} has $\pi^\star\equiv0$ on $\{R<\beta^\star\}$. Suppose further that $\mu$ puts mass at least $\omega>0$ on reduced states $R\ge\xi$, for some $\xi\in\cS$ with $\xi>\varepsilon$, and that
\begin{equation}
  \limsup_{R\downarrow\varepsilon}\ \Prob\bigl[R_1>\varepsilon \mid R_0=R\bigr]
  \;<\; \frac{\omega\,T^\star(\xi)}{\bar T},
  \label{eq:plow_condition}
\end{equation}
where $\bar T$ is a uniform upper bound on $T^\star$,
such as~\eqref{eq:Tbar} when applicable, and
$R_1\sim P(\cdot\mid R)$. Then $\beta^\star>\varepsilon$:
the no-restart region is non-degenerate.
\end{corollary}

The proof is in Appendix~\ref{app:proof_threshold}; Remark~\ref{rem:plow_needed} explains why a condition such as~\eqref{eq:plow_condition} cannot be dropped, and Remarks~\ref{rem:calib} and~\ref{rem:randomise} discuss the calibration of $\beta^\star$ and the case for randomisation.

\section{A state-dependent stochastic restart rule}
\label{sec:policy}

The rule evaluated in Section~\ref{sec:experiments}, which we call sPDLP \footnote{The implementation is available at
\url{https://github.com/ykoubaa-code/sPDLP.git}.}, is stated in Algorithm~\ref{alg:stochastic_restart}. It is evaluated every $40$ iterations, at the termination check of PDLP, and it replaces the adaptive restart of PDLP; the step-size and primal-weight updates are unchanged. 
It uses the running-average score
$\bar s^k:=s(\bar z^k)$, with $s$ defined in ~\eqref{eq:score},
to suppress restarting near convergence, and the \emph{merit}
$M_k 
=\sqrt{R_{k-1}}$,
where $R_k$ is defined in Section~\ref{sec:setting}, to modulate the probability away from it, together with the \emph{epoch length} $L_k$, the number of iterations since the last restart, and the \emph{restart count} $\nu_k$, the number of stochastic restarts fired so far. Motivated by Corollary~\ref{cor:plow}, the rule sets the restart probability to zero in the region

\begin{equation}
    \bar s^k \le \kappa\,\varepsilon_{\mathrm{opt}},
    \qquad \kappa=2.5,
    \label{eq:no_restart_region}
\end{equation}
where $\varepsilon_{\mathrm{opt}}$ is the termination tolerance. Outside this region, the probability is
\begin{equation}
p_k =
\max\bigl\{1-\exp(-\eta\,\zeta_k),\,p_{\min}\bigr\},
\qquad \eta = -\log(1-p_{\max}),
\label{eq:restart_probability}
\end{equation}
with the activation $\zeta_k$ built from the merit, the epoch length and the restart count:
\begin{align}
\chi_k &= 1-e^{-M_k}, &
L_k^{\mathrm{eff}} &= L_0\sqrt{M_k}\;2^{\,0.6\,\lfloor \nu_k/d\rfloor}, \label{eq:reff}\\
\alpha_k &= \frac{L_k}{L_k+L_k^{\mathrm{eff}}}, &
\zeta_k &= \alpha_k\bigl[\chi_k+(1-\chi_k)\alpha_k\bigr], \label{eq:activation}
\end{align}
where $L_0=200$ and $d=3$. Since $\alpha_k\in[0,1]$, one has
$\alpha_k^2\le\zeta_k\le\alpha_k
<1$
and so $p_k\in[p_{\min},p_{\max}]$. The patience factor $\alpha_k$ allows each epoch a window of order $L_k^{\mathrm{eff}}$ before a restart becomes likely, and the window lengthens by about one half for every $d$ restarts taken. The no-restart condition ~\eqref{eq:no_restart_region} is motivated by Corollary~\ref{cor:plow} but uses the running-average convergence score $\bar s^k$ rather than the theoretical state $R$, to align with PDLP's stopping criterion. This score-based condition is heuristic and is not covered by the optimality guarantee of Section~\ref{sec:threshold}. The remaining rule is also heuristic: its dependence on $L_k$ and $\nu_k$ lies outside the reduced-state
control model, while the floor $p_{\min}>0$ departs
from the zero restart probability of Corollary~\ref{cor:plow}.

\begin{algorithm}[H]
\caption{sPDLP: the state-dependent stochastic restart rule, evaluated every $40$ iterations}
\label{alg:stochastic_restart}
\begin{algorithmic}[1]
\Require iterates $z^k,z^{k-1}$; average $\bar z^k$; restart count $\nu_k$; epoch length $L_k$; parameters $\kappa,d,L_0,p_{\max},p_{\min}$
\State $\bar s^k \gets s(\bar z^k)$ by \eqref{eq:score}
\If{$\bar s^k \le \kappa\,\varepsilon_{\mathrm{opt}}$} \Return no restart \Comment{no-restart region}
\EndIf
\State $M_k \gets \left(
\frac{1}{\tau_{k-1}}\|x^k-x^{k-1}\|_2^2
+
\frac{1}{\sigma_{k-1}}\|y^k-y^{k-1}\|_2^2
\right)^{1/2}
$;\quad compute $p_k$ from \eqref{eq:restart_probability}--\eqref{eq:activation}
\State draw $v\sim\mathrm{Unif}(0,1)$;\quad \textbf{if} $v>p_k$ \textbf{then} \Return no restart
\State draw $\theta\sim\mathrm{Unif}(0,1)$;\quad $z^{\mathrm{dest}} \gets
\operatorname*{arg\,min}_{
z\in\{z^k,\bar z^k,
\theta\bar z^k+(1-\theta)z^k\}} s(z)$
\State move to $z^{\mathrm{dest}}$; reset the running average; $\nu_{k+1}\gets \nu_k+1$, $L_{k+1}\gets0$
\end{algorithmic}
\end{algorithm}

When a restart fires, the destination is the lowest-scoring of the current iterate, the running average and a random interpolation between them,
\begin{equation}
z^{\mathrm{dest}}
=
\arg\min\bigl\{s(z):\ z\in\{z^k,\;\bar z^k,\;\theta\bar z^k+(1-\theta)z^k\}\bigr\},
\qquad \theta\sim\mathrm{Unif}(0,1).
\label{eq:best_destination}
\end{equation}
Since $z^k$ is itself a candidate, a restart never increases the score; when $z^k$ is chosen only the averaging window is reset. The residuals of the interpolated candidate follow by linearity from $Kx^k$, $K\bar x^k$, $K^\top y^k$ and $K^\top\bar y^k$, so the rule needs no product with $K$ beyond those of the termination check.

\paragraph{Snapping.}
Every destination in~\eqref{eq:best_destination} lies in
the convex hull of the trajectory and
cannot bypass the plateau of
Theorem~\ref{thm:wall}(i). To overcome this limitation, when a restart
fires, we also consider $(\hat x,y^k)$ in the solver's scaled space,
where each coordinate of $x^k$ is independently moved with
probability one half to its upper bound if
$(c-K^\top y^k)_j<0$ and $u_j<\infty$, or to its lower bound if
$(c-K^\top y^k)_j>0$ and $l_j>-\infty$, and left unchanged otherwise.
The candidate joins the minimisation
in~\eqref{eq:best_destination}, so a restart still
never increases the score;
coordinates already at the indicated
bound remain unchanged.
The extra cost is one product with $K$ per restart
that fires.
We call the result sPDLP with snapping.
Figures~\ref{fig:pagerank-iterations} and~\ref{fig:wall}
use it; the main benchmark of
Section~\ref{sec:experiments} uses
Algorithm~\ref{alg:stochastic_restart} without snapping.
Results with snapping on the same benchmark appear
in Appendix~\ref{app:snapping_benchmark}.

\section{Computational results}
\label{sec:experiments}

\paragraph{Setup.}
We use a benchmark set of $68$ instances from the Mittelmann collection, LP relaxations of MIPLIB 2017 models \citep{gleixner2021miplib}, whose sizes after presolve are listed in Table~\ref{tab:main}. The baseline is PDLP as implemented in the \texttt{FirstOrderLp.jl} code released by Google Research, with its default adaptive normalised-duality-gap restart; sPDLP is the same code with that restart replaced by Algorithm~\ref{alg:stochastic_restart}. Each run targets $\varepsilon_{\mathrm{opt}}=10^{-4}$ and is capped at $3600$ seconds or $10^6$ iterations; sPDLP is run with three seeds per instance. The iteration ratio of an instance, and likewise its time ratio, is the geometric mean over seeds of the sPDLP count divided by the PDLP count, so that a ratio below one favours sPDLP; the overall figure is the geometric mean over instances. On \texttt{map18} sPDLP reaches the time limit on all three seeds while PDLP terminates; on \texttt{rail02} PDLP reaches the time limit while sPDLP terminates on one seed; ratios are therefore reported on the remaining $66$ instances. The hyperparameters, $\kappa=2.5$, $d=3$, $L_0=200$, the exponent $0.6$ in~\eqref{eq:reff}, $p_{\max}=0.7$ and $p_{\min}=0.003$, were chosen by a grid search on fourteen of the $68$ instances, listed in Appendix~\ref{app:table}, and then fixed for all.

\paragraph{Results.}
Table~\ref{tab:main} in the appendix gives the full results and Table~\ref{tab:summary} their distribution; Figure~\ref{fig:gmratio} shows the per-instance iteration ratios with their variability across seeds. In iterations the rule is close to neutral: the geometric-mean ratio is $1.074$, sPDLP needs fewer iterations on $27$ of the $66$ instances and more on $39$, and $58$ of the $66$ ratios lie within a factor of two of one. A two-sided sign test on the $27$ wins gives $p\approx0.18$, so the benchmark does not distinguish the two methods. The one large effect is \texttt{app1-2}, the instance on which the baseline is by far the slowest, where sPDLP needs $4.5$ per cent of the baseline's iterations; without it the geometric mean is $1.128$. On the losing side, \texttt{glass4} is the only pronounced outlier, with a ratio of $4.421$; the remaining large losses are factors between $2.0$ and $2.7$. In wall-clock time the geometric-mean ratio is $1.078$ and sPDLP is faster on $33$ of the $66$ instances. On several instances the time ratio and the iteration
ratio differ widely, potentially reflecting compilation overheads and runtime variability. The wall-clock comparisons should therefore be interpreted with caution.

\begin{figure}[H]
\centering
\includegraphics[width=\linewidth]{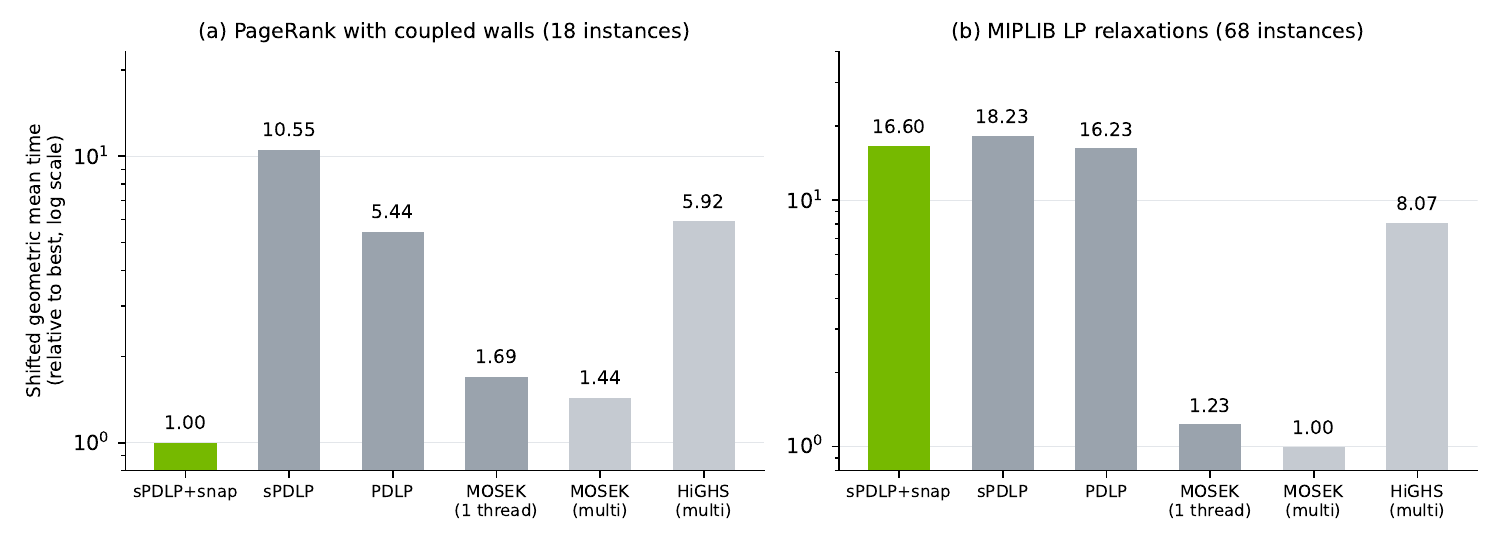}
\caption{Shifted geometric-mean solve times
(10\,s shift, log scale), normalised to the best solver
per panel: (a) 18 coupled-wall PageRank LPs
(Appendix~\ref{app:pagerank});
(b) 68 MIPLIB LP relaxations
(Section~\ref{sec:experiments}).
Timeouts, out-of-memory runs, objective errors above $10^{-3}$, and larger instances not attempted after failures on smaller ones are penalised at $3600$\,s.
For unsuccessful runs, the time and iteration counts are set to their respective limits before aggregation
across seeds.
sPDLP with snapping uses three seeds; sPDLP uses
one in (a) and three in (b).
Dark and light bars indicate single- and multithreaded
runs. Solved counts (left to right):
(a) $18,0,5,9,9,3$; (b) $68,66,67,68,68,67$.}
\label{fig:solver-comparison}
\end{figure}

\paragraph{Coupled-wall PageRank experiments.}
The PageRank experiments of
Appendix~\ref{app:pagerank} illustrate a regime
favourable to snapping. With fixed steps and no
diagonal preconditioning, sPDLP with snapping
requires up to $13.5$ times fewer iterations than
PDLP on completed one-wall comparisons.
Figure~\ref{fig:solver-comparison} compares
penalised aggregate solve times on PageRank and
MIPLIB: snapping performs best on PageRank,
whereas multithreaded MOSEK performs best
on MIPLIB under this measure.

\paragraph{Relation to the theory.}
Section~\ref{sec:exact} concerns geometric restarts
superimposed on a fixed base method, whereas sPDLP
replaces the restart of PDLP by a state-dependent rule,
so the benchmark is not a direct test of the theory.
The large gain on \texttt{app1-2}
is consistent with the slow-mode intuition of
Remark~\ref{rem:notnecessary}, but does not establish
the assumptions of the renewal model.
On \texttt{danoint}, a controlled test with PDHG
as the base method, the renewal formula
predicts a slowdown under constant-rate
restarting, consistent with the observations,
but overestimates its magnitude
(Appendix~\ref{app:danoint}).
The 200 restart-free runs yield an
empirical NBUE distribution, consistent with
the observed slowdown under geometric restarting
(Corollary~\ref{cor:nbue} and
Appendix~\ref{app:danoint-nbue}).

\section{Conclusion}

Randomised restarting can beat PDLP by a margin that grows with the dimension, as the plateau-and-wall family shows, but only when a restart can leave the convex hull of the trajectory, as snapping does. When a run that has not yet finished can expect to finish no later than a fresh run would, no constant restart rate helps; when the iteration count varies little from one restart point to another, as under a two-sided linear rate, none helps appreciably; and in a reduced-state model with monotone
dynamics the optimal policy
has a two-level threshold structure, with a
non-degenerate no-restart neighbourhood under
an additional condition.
A natural next step is a restart rule closer to the theory, which would choose its probability online from estimates of $T_p$ (Section~\ref{sec:practice}) and allow for the dependence of the restart destination on the trajectory.

\section*{Acknowledgements}
This work has been supported by European Union's HORIZON–MSCA-2022-DN-JD programme under the Horizon Europe (HORIZON) Marie Skłodowska-Curie Actions, grant agreement 101120296 (TENORS). 

\section*{Reproducibility statement}
Proofs of all results are given in Appendices~\ref{app:proofs}, \ref{app:proof_threshold} and~\ref{app:proof_wall}, with the diffusion surrogate in Appendix~\ref{app:surrogate}. Algorithm~\ref{alg:stochastic_restart} states the rule in full, and Section~\ref{sec:experiments} gives every hyperparameter, the tolerance, the limits and the number of seeds. Appendix~\ref{app:table} lists the tuning instances, the hardware and the per-instance results. The instances of Figure~\ref{fig:wall} are fully specified in Section~\ref{sec:wall}, with the family parameters in the caption, and the PageRank instances in Appendix~\ref{app:pagerank}. The baseline is the publicly released \texttt{FirstOrderLp.jl} code. %

\section*{The use of large language models}
A large language model was used to revise the wording and structure of the manuscript, to check the consistency of its notation, to help search the literature on first passage under restart, and to help formulate the example of Theorem~\ref{thm:wall}. The authors have checked all text and all proofs and take full responsibility for the content of the paper. %

\bibliographystyle{plainnat}
\bibliography{Reference}

\newpage
\appendix
\raggedbottom
\section{Idealised plateau-and-wall constructions}
\label{app:plateau-experiments}

\subsection{Numerical illustration and limitations}
\begin{figure}[H]
\centering
\includegraphics[width=0.92\textwidth]{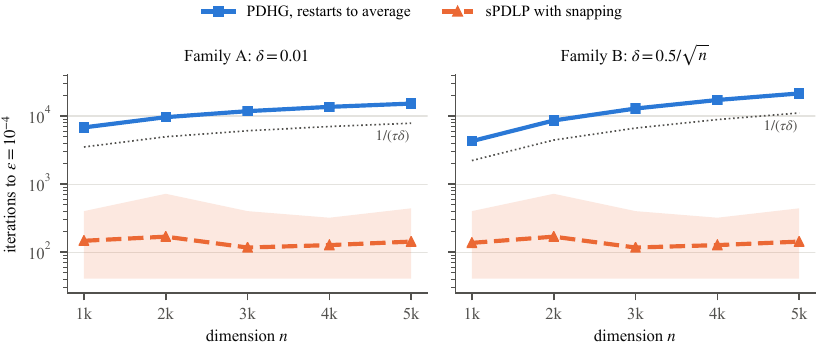}
\caption{Iterations to reach $\varepsilon=10^{-4}$ on plateau-and-wall instances with $n$ from $1000$ to $5000$, in a reference implementation of PDHG with fixed step sizes and no preconditioning. Restarts to the running average, the destination PDLP uses, here every $64$ iterations, roughly double the plateau length $1/(\tau\delta)$ (dotted). sPDLP with snapping (Section~\ref{sec:policy}) needs between $117$ and $169$ iterations on average over $20$ seeds at every $n$; the band shows the range over seeds. 
In family B, $\delta$ shrinks like
$1/\sqrt{n}$, making the wall gradient smaller
as the dimension increases.
}
\label{fig:wall}
\end{figure}

To overcome the plateau, sPDLP with snapping adds a
fourth restart candidate that moves selected variables
directly to their bounds, guided by their reduced
costs (Section~~\ref{sec:policy}). Figure~\ref{fig:wall}
shows that this variant needs about 150 iterations
across dimensions, whereas PDHG with averaging
restarts needs thousands, growing with $n$.
The example is idealised: presolve would eliminate
the unconstrained wall variable, diagonal
preconditioning can reduce the factor $\sqrt{n}$,
and the plateau matters only while $\delta$ is
large relative to the tolerance. The
dimension-dependent plateau also lies outside
the bounded-spread regime discussed after
Corollary~\ref{cor:linear}. Appendix~\ref{app:ml} discusses connections
to gradient behaviour in ReLU networks.

\subsection{Multiple walls and coordinate-wise restarting}
\label{app:many-walls}

We extend the construction of
Section~1.1 to multiple walls and compare
global with coordinate-wise restarting.

Take $n_w$ variables $w=(w_1,\dots,w_{n_w})\in[0,1]^{n_w}$, each with cost $-\delta$ and entering no constraint, together with the block of Section~\ref{sec:wall}, of norm $\|\hat K\|_2=\sqrt n$, and run PDHG with $\tau=\sigma$ and $\tau\|\hat K\|_2<1$. Since $w$ is decoupled, each $w_i$ moves to $\min\{1,w_i+\tau\delta\}$ at every step, independently of everything else. The target is $w=\mathbf 1$ together with a target set $\mathcal B$ for the block, and we assume that from every starting point the block iterates enter $\mathcal B$ within $f$ iterations and never leave it. The second half of this assumption holds when $\mathcal B$ is a ball about the solution set in the norm in which PDHG is Fejér monotone \citep{chambolle2011first}; the first half is the idealisation. Call $w_i$ \emph{settled} when $w_i=1$: it sits at a bound, and its reduced cost $-\delta$ holds it there, a test that uses only the sign of the gradient. We compare three schemes, with $L:=\lceil1/(\tau\delta)\rceil$ as before.
\begin{enumerate}
\item[(T)] Methods whose iterates are PDHG steps or convex combinations of earlier iterates, started at $w=\mathbf 0$.
\item[(G)] Global geometric restarting at rate $p$: with probability $p$ at each iteration the whole iterate is replaced by a fresh draw, with $w_1,\dots,w_{n_w}$ independent and uniform on $\{0,1\}$ and the block from any fixed law.
\item[(C)] Coordinate-wise restarting at rate $p\in(0,1]$: after each PDHG step, each unsettled $w_i$ is replaced, independently and with probability $p$, by a uniform draw from $\{0,1\}$; the block is never reset.
\end{enumerate}

\begin{proposition}[Many walls]\label{prop:manywalls}
Let $n_w\ge1$ and $f\ge1$.
\begin{enumerate}
\item[(i)] Every method of type (T) needs at least $L\ge\sqrt n/\delta$ iterations.
\item[(ii)] Under (G), for every $p\in(0,1)$, the expected number of iterations is at least $\tfrac14\min\{L,\,2^{n_w}\}$.
\item[(iii)] Under (C), from any starting point, the expected number of iterations is at most $f+1+\tfrac2p\,(1+\ln n_w)$.
\end{enumerate}
\end{proposition}

\begin{proof}
(i) is Theorem~\ref{thm:wall}(i) applied to $w_1$.

(ii) Fresh draws make the restart-free counts of successive epochs independent, so Assumption~\ref{ass:N} holds and, by the proof of Theorem~\ref{thm:renewal}, $T_p=\E[\min(N,\Gamma)]/\Prob[\Gamma\ge N]$. Let $\varpi:=2^{-n_w}\le\tfrac12$ and $u:=1-p$. With probability $\varpi$ the draw has $w=\mathbf 1$; otherwise some $w_i$ starts at $0$, and since it rises by exactly $\tau\delta$ per step it needs $L$ steps, so $N\ge L$. Hence $\Prob[\Gamma\ge N]\le\varpi+u^{L-1}$ and $\E[\min(N,\Gamma)]\ge\max\bigl\{1,\ \tfrac12\E[\min(L,\Gamma)]\bigr\}$, where $\E[\min(L,\Gamma)]=\sum_{t=0}^{L-1}u^t$. If $u^{L-1}\le\varpi$, then $T_p\ge1/(2\varpi)=2^{n_w-1}$. Otherwise $\Prob[\Gamma\ge N]\le2u^{L-1}$ and
\[
T_p\;\ge\;\frac{\sum_{t=0}^{L-1}u^t}{4\,u^{L-1}}\;=\;\frac14\sum_{i=0}^{L-1}u^{-i}\;\ge\;\frac L4 .
\]

(iii) Let $\Theta_i$ be the first iteration at which $w_i=1$. A settled coordinate stays settled, since the PDHG step keeps it at $1$ and it is never reset. Before $\Theta_i$, at every iteration, conditionally on the past, $w_i$ is reset with probability $p$ and then set to $1$ with probability $\tfrac12$, so $\Prob[\Theta_i>t]\le(1-p/2)^t$. By the union bound, $\Prob[\max_i\Theta_i>t]\le n_w(1-p/2)^t$. With $t_0:=\lceil(2/p)\ln n_w\rceil$, so that $n_w(1-p/2)^{t_0}\le n_w\mathrm e^{-pt_0/2}\le1$,
\[
\E\bigl[\max_i\Theta_i\bigr]=\sum_{t\ge0}\Prob\bigl[\max_i\Theta_i>t\bigr]\le t_0+\sum_{t\ge t_0}n_w(1-p/2)^t\le\frac2p\ln n_w+1+\frac2p .
\]
The block is never reset, so it follows plain PDHG, enters $\mathcal B$ within $f$ iterations and stays there. The target is therefore reached by iteration $\max\{f,\max_i\Theta_i\}\le f+\max_i\Theta_i$.
\end{proof}

\begin{table}[h]
\caption{Iteration counts on the many-walls model of Proposition~\ref{prop:manywalls}, with $n_w$ walls, a block of norm $\sqrt n$ that settles within $f$ iterations, and $p=\tfrac12$ in scheme (C).}
\label{tab:manywalls}
\centering
\footnotesize
\begin{tabular}{lll}
\toprule
Scheme & What a restart replaces & Expected iterations \\
\midrule
(T) PDLP-type restarts & the iterate, by a convex combination of past iterates & at least $\sqrt n/\delta$ \\
(G) global restarts, any rate & the whole iterate, by a fresh draw & at least $\tfrac14\min\{\sqrt n/\delta,\ 2^{n_w}\}$ \\
(C) coordinate-wise restarts & each unsettled wall coordinate & at most $f+5+4\ln n_w$ \\
\bottomrule
\end{tabular}
\end{table}

Table~\ref{tab:manywalls} collects the bounds. With $n_w$ walls and $2^{n_w}\le\sqrt n/\delta$, restarting the whole iterate is slower than restarting stuck coordinates by a factor of order $2^{n_w}/(f+\ln n_w)$. Resetting a whole iterate at random is the mechanism of the stochastic-resetting literature, but a restart that lands on all the walls at once becomes exponentially unlikely as their number grows, so the whole iterate is the wrong unit to reset.

\clearpage
\section{PageRank LPs with Coupled Walls}\label{app:pagerank}

\paragraph{Motivation.}
To see whether the benefit of snapping survives when the walls are coupled to the rest of the problem, we place them in a PageRank LP on a graph with two communities. Each wall moves a capped amount of the random surfer's restart probability from a page in one community to a page in the other, and so changes both the PageRank equations and the PageRank vector. Interventions of this kind have been studied on networks that fall into two communities with few links between them, such as the two sides of a polarised debate or pages written in two languages. \citet{tsioutsiouliklis2021fairness} alter the jump vector of PageRank so that a protected group receives a prescribed share of the PageRank, and \citet{garimella2017reducing} add links across the divide of a polarised network so as to reduce a random-walk measure of controversy. Our objective is deliberately simpler than theirs, so that the optimum is known in closed form. The family tests the mechanism of Theorem~\ref{thm:wall}; we do not claim the corresponding lower bound for the coupled problem.

We use the PageRank LP of \citet{nesterov2014subgradient}, normalised by $\mathbf 1^\top x=1$, on Barab\'asi--Albert graphs of the kind with which \citet{Applegate2021PracticalLL} demonstrated the scalability of PDLP. The normalisation row has Euclidean norm $\sqrt N$, so that $\|K\|_2\ge\sqrt N$. Without preconditioning the fixed step is $0.8/\hat\sigma$, where $\hat\sigma$ estimates $\|K\|_2$ by power iteration, and the plateau of Theorem~\ref{thm:wall} would then last about $1/(\tau\delta)\approx1.25\sqrt N/\delta$ iterations. For the coupled problem this is a heuristic scale, not a bound.

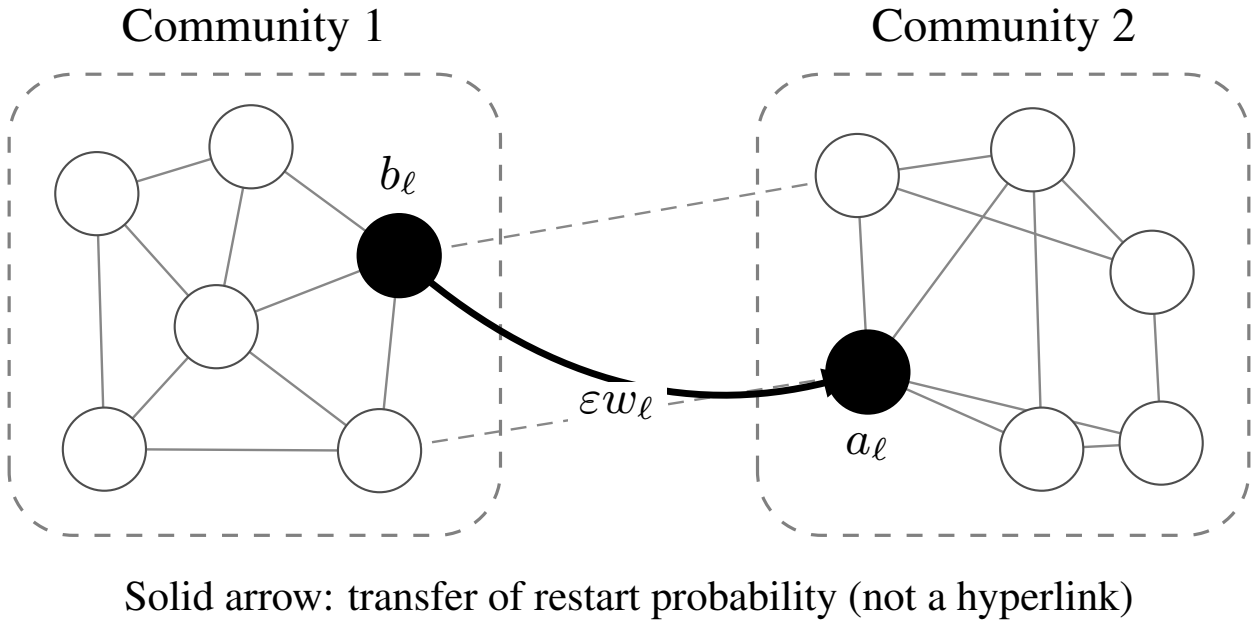
\begin{figure}[H]
\centering
\begin{tikzpicture}[
    >=Latex,
    ordinary/.style={line width=0.45pt,draw=black!47},
    crosslink/.style={ordinary,densely dashed},
    page/.style={circle,draw=black!70,fill=white,minimum size=5.5mm,inner sep=0pt},
    endpoint/.style={circle,draw=black,fill=black,minimum size=5.5mm,inner sep=0pt},
    cluster/.style={draw=black!50,dashed,rounded corners=12pt,line width=0.6pt}]
  \draw[cluster] (-4.1,-1.50) rectangle (-0.85,1.55);
  \draw[cluster] (0.85,-1.50) rectangle (4.1,1.55);
  \coordinate (u1) at (-3.52, 0.76);
  \coordinate (u2) at (-2.50, 1.07);
  \coordinate (b)  at (-1.52, 0.35);
  \coordinate (u3) at (-2.73,-0.12);
  \coordinate (u4) at (-3.47,-0.93);
  \coordinate (u5) at (-1.65,-0.94);
  \coordinate (v1) at (1.51,0.88);
  \coordinate (v2) at (2.67,1.05);
  \coordinate (v3) at (3.46,0.24);
  \coordinate (a)  at (1.58,-0.42);
  \coordinate (v4) at (2.73,-0.93);
  \coordinate (v5) at (3.52,-0.89);
  \foreach \s/\t in {u1/u2,u2/b,b/u3,u3/u4,u4/u5,
                        u5/b,u1/u3,u2/u3,u3/u5,u1/u4,
                        v1/v2,v2/v3,v3/v5,v5/v4,
                        v4/a,a/v1,v1/v3,v2/a,a/v5,v2/v4}
    \draw[ordinary] (\s) -- (\t);
  \draw[crosslink] (b) -- (v1);
  \draw[crosslink] (u5) -- (a);
  \foreach \p in {u1,u2,u3,u4,u5,v1,v2,v3,v4,v5}
    \node[page] at (\p) {};
  \node[endpoint,label={[font=\small]above:$b_\ell$}] at (b) {};
  \node[endpoint,label={[font=\small]below:$a_\ell$}] at (a) {};
  \node[font=\small] at (-2.48,1.84) {Community 1};
  \node[font=\small] at ( 2.48,1.84) {Community 2};
  \draw[-{Latex[length=3mm,width=2mm]},line width=1.3pt]
      (b) to[bend right=28]
      node[below,fill=white,inner sep=2pt,font=\small]
      {$\varepsilon w_\ell$} (a);
  \node[font=\footnotesize,align=center] at (0,-1.93)
      {Solid arrow: transfer of restart probability (not a hyperlink)};
\end{tikzpicture}
\caption{Two-community PageRank with coupled walls (schematic, not a drawing
of a generated instance). The communities have many internal graph links
and few cross-community links (dashed). A wall transfers
$\varepsilon w_\ell$ units of \emph{restart} probability from $b_\ell$
to $a_\ell$. It changes the PageRank equations without adding an edge
to the hyperlink graph.}
\label{fig:pagerank-communities}
\end{figure}

\paragraph{Construction.}
Let $G$ comprise two Barab\'asi--Albert graphs of $N/2$ nodes each,
with attachment parameter $3$, joined by $50$ randomly chosen
cross-community edges. Let $S\in\mathbb R^{N\times N}$ be the
column-stochastic random-walk matrix with
$S_{ij}=1/\deg(j)$ when $\{i,j\}$ is an edge, and zero otherwise;
let $d=0.99$. For each wall $\ell=1,\ldots,k$, choose distinct
endpoints $b_\ell$ in the first community and $a_\ell$ in the second,
and set
\[
 B=\bigl[e_{b_1}-e_{a_1},\ldots,e_{b_k}-e_{a_k}\bigr]
 \in\mathbb R^{N\times k}.
\]
The coupled-wall LP is
\begin{equation}\label{eq:pagerank-coupled-walls}
\begin{aligned}
\min_{x\in\mathbb R^N,\,w\in\mathbb R^k}\quad
 &-\delta\sum_{\ell=1}^k w_\ell \\
\text{subject to}\quad
 &(I-dS)x+\varepsilon Bw\geq\frac{1-d}{N}\mathbf{1},\\
 &\mathbf{1}^{\top}x=1,\\
 &0\leq w\leq\mathbf{1},\qquad x\text{ free},
\end{aligned}
\end{equation}
where $\delta=10^{-3}$ and
$\varepsilon=\rho(1-d)/N$ with $\rho=\tfrac12$.
Each wall modifies two PageRank equations; increasing $w_\ell$
transfers $\varepsilon w_\ell$ units of teleportation probability
from $b_\ell$ to $a_\ell$, without adding an edge to the graph.
Figure~\ref{fig:pagerank-communities} illustrates the construction.

\paragraph{Feasibility and known optimum.}
Since $\mathbf 1^\top S=\mathbf 1^\top$ and
$\mathbf 1^\top B=0$, summing the PageRank inequalities and using
$\mathbf 1^\top x=1$ forces all inequalities to hold at equality.
For any $w\in[0,1]^k$, this yields the unique PageRank vector
\begin{equation}\label{eq:pagerank-x-of-w}
 x(w)=(I-dS)^{-1}\left(\frac{1-d}{N}\mathbf 1-\varepsilon Bw\right).
\end{equation}
The inverse $(I-dS)^{-1}=\sum_{t=0}^\infty d^tS^t$
is entrywise nonnegative. Because all $2k$ wall endpoints are
distinct and $\rho=\tfrac12$, the vector in parentheses in
\eqref{eq:pagerank-x-of-w} is componentwise at least
$(1-\rho)(1-d)\mathbf 1/N>0$. Thus $x(w)>0$ for
all $w\in[0,1]^k$ and every wall setting is feasible.
The unique optimum is
\[
 w^\star=\mathbf 1,\qquad x^\star=x(\mathbf 1),
 \qquad f^\star=-k\delta.
\]
The LP has $N+1$ constraint rows, $N+k$ columns, and
approximately $8N$ nonzero coefficients.

\paragraph{Wall reduced costs.}
At any dual-feasible point $(y,\mu)$, since $x$ is free
and has zero objective coefficient,
$(I-dS^\top)y+\mu\mathbf 1=0$. Because
$S^\top\mathbf 1=\mathbf 1$ and $I-dS^\top$ is nonsingular,
$y=-\mu\mathbf 1/(1-d)$. It follows that the wall reduced costs are
\begin{equation}\label{eq:pagerank-reduced-cost}
 r_w=-\delta\mathbf 1-\varepsilon B^\top y
     =-\delta\mathbf 1.
\end{equation}
The coupling thus alters primal feasibility but leaves the reduced costs of the walls unchanged at every \emph{dual-feasible} point; the iterates of PDHG, of course, need not be dual feasible.

\paragraph{Experimental configuration.}
We compare PDLP, sPDLP and sPDLP with snapping in the setting of Theorem~\ref{thm:wall}, not in the default configuration of PDLP: Ruiz, Pock--Chambolle and $\ell_2$ rescaling are switched off, the step size is fixed, and the primal weight is fixed at $1$. The stopping tolerances are $10^{-4}$, relative and absolute, with limits of $3600$\,s and $10^6$ iterations. The setting is chosen to expose the plateau; diagonal preconditioning rescales the wall columns and can reduce the effect.

MOSEK uses its interior-point method
with presolve, and HiGHS uses dual simplex;
other solver settings are left at their defaults,
apart from the stated tolerance.

The PDLP variants run on one thread, in Julia, and each experiment was run on its own on a laptop with an Intel Core i5-1145G7 and 16\,GB of memory. The times of the PDLP variants include the compilation of Julia code but not the reading of the MPS file; those of MOSEK and HiGHS are as reported by the optimisers, so the two are not measured in quite the same way. The snapping rule is described in Section~\ref{sec:policy}.
The snapping results are geometric means over seeds $11$, $22$, and $33$;
all $54$ reported snapping runs end with $w=\mathbf 1$.

\paragraph{Results.}
Table~\ref{tab:c4-coupled-walls} gives the full results, and Figure~\ref{fig:pagerank-iterations} plots the one-wall comparisons that terminate. At $N=10^5$ and $k=1$, sPDLP with snapping takes $28{,}591$ iterations, against $385{,}000$ for PDLP and $415{,}400$ for sPDLP without snapping, factors of about $13.5$ and $14.5$. At $N=3\cdot10^5$ and $k=1$, both variants without snapping reach the time limit, whereas the variant with snapping terminates. At $N=10^5$, sPDLP with snapping reports $213.9$\,s against $740.4$\,s for multithreaded MOSEK, although the two times are not measured in the same way.

\paragraph{Why runs stop short of the wall.}
At a feasible $w$ the objective gap is
\begin{equation}\label{eq:pagerank-objective-gap}
 f(w)-f^\star=\delta\sum_{\ell=1}^k(1-w_\ell).
\end{equation}
On the plateau the primal and dual residuals are already small, and since the objective values are of order $k\delta\le0.032$, the gap term of the stopping test is close to~\eqref{eq:pagerank-objective-gap}. A run that has not reached the wall therefore stops once $\delta\sum_\ell(1-w_\ell)$ has fallen to about $10^{-4}$. The walls start together and move at nearly the same speed, so this happens at $w_\ell\approx1-10^{-4}/(k\delta)=1-0.1/k$, that is, at $0.9$, $0.9875$ and $0.997$ for $k=1$, $8$ and~$32$; the bracketed values in Table~\ref{tab:c4-coupled-walls} are at or just above these levels. Such a run meets the stopping test, yet its objective error relative to $|f^\star|=k\delta$ is about $0.1/k$, which exceeds the threshold of $10^{-3}$ used in Figure~\ref{fig:solver-comparison} for every $k$ in the table; it is therefore counted there as a failure. Every run with snapping ends at $w^\star=\mathbf 1$.

\begin{table}[H]
\centering
\caption{
Two-community PageRank LPs with $k$ coupled walls ($\delta=10^{-3}$).
Iteration counts and solve times (seconds). sPDLP+snap reports geometric
means over three random seeds. $^\dagger$ indicates termination with some
$w_\ell<1$ (the bracketed number gives the smallest $w_\ell$); at a time limit,
brackets give the smallest value reached. Bold marks the fastest completed
run at each size. ``t.l.'': 3600\,s limit; ``o.o.m.'': out of memory
(16\,GB); ``--'': not run. The PDLP variants use fixed steps and no
preconditioning; MOSEK uses interior-point and
HiGHS dual simplex; other settings are left
at their defaults.
Missing comparisons are not treated as completed runs.}
\label{tab:c4-coupled-walls}
\small
\setlength{\tabcolsep}{4pt}
\begin{tabular}{rr rr rr rr rr r}
\toprule
 & & \multicolumn{2}{c}{sPDLP+snap} & \multicolumn{2}{c}{sPDLP} & \multicolumn{2}{c}{PDLP}
 & \multicolumn{2}{c}{MOSEK} & HiGHS \\
\cmidrule(lr){3-4}\cmidrule(lr){5-6}\cmidrule(lr){7-8}\cmidrule(lr){9-10}\cmidrule(lr){11-11}
$N$ & $k$ & it. & time & it. & time & it. & time & 1 thr. & multi & time \\
\midrule
$10^4$ & 1 & 14{,}560 & 13.7 & 150{,}560$^\dagger$ {\scriptsize[0.900]} & 91.3 & 118{,}120$^\dagger$ {\scriptsize[0.929]} & 107.1 & 2.1 & \textbf{1.3} & 105.4 \\
 & 8 & 14{,}328 & 14.9 & 131{,}440$^\dagger$ {\scriptsize[0.987]} & 118.4 & 128{,}720 & 113.9 & 2.9 & \textbf{2.3} & 108.4 \\
 & 32 & 14{,}870 & 14.4 & 128{,}800$^\dagger$ {\scriptsize[0.997]} & 89.4 & 128{,}520 & 101.0 & 2.3 & \textbf{1.3} & 98.7 \\
\midrule
$3\cdot 10^4$ & 1 & 21{,}742 & 45.2 & 236{,}280$^\dagger$ {\scriptsize[0.900]} & 561.0 & 223{,}680 & 584.8 & 43.2 & \textbf{30.6} & t.l. \\
 & 8 & 20{,}931 & 60.0 & 226{,}880$^\dagger$ {\scriptsize[0.987]} & 483.1 & 221{,}120 & 707.2 & 32.1 & \textbf{24.4} & t.l. \\
 & 32 & 21{,}872 & 67.1 & 222{,}040$^\dagger$ {\scriptsize[0.997]} & 583.6 & 221{,}120 & 716.2 & 44.3 & \textbf{31.4} & t.l. \\
\midrule
$10^5$ & 1 & 28{,}591 & \textbf{213.9} & 415{,}400$^\dagger$ {\scriptsize[0.900]} & 2677.3 & 385{,}000$^\dagger$ {\scriptsize[0.900]} & 3205.6 & 1416.8 & 740.4 & -- \\
 & 8 & 28{,}195 & \textbf{208.4} & -- & -- & -- & -- & 1371.1 & 705.5 & -- \\
 & 32 & 29{,}693 & \textbf{223.1} & -- & -- & -- & -- & 1428.8 & 740.6 & -- \\
\midrule
$3\cdot 10^5$ & 1 & 35{,}384 & \textbf{975.0} & 139{,}040 {\scriptsize[0.143]} & t.l. & 113{,}680 {\scriptsize[0.127]} & t.l. & -- & o.o.m. & -- \\
 & 8 & 34{,}276 & \textbf{935.1} & -- & -- & -- & -- & -- & -- & -- \\
 & 32 & 35{,}305 & \textbf{965.6} & -- & -- & -- & -- & -- & -- & -- \\
\midrule
$5\cdot 10^5$ & 1 & 40{,}566 & \textbf{1881.5} & -- & -- & -- & -- & -- & -- & -- \\
 & 8 & 39{,}236 & \textbf{1812.9} & -- & -- & -- & -- & -- & -- & -- \\
 & 32 & 40{,}400 & \textbf{1880.8} & -- & -- & -- & -- & -- & -- & -- \\
\midrule
$7\cdot 10^5$ & 1 & 41{,}330 & \textbf{2807.7} & -- & -- & -- & -- & -- & -- & -- \\
 & 8 & 41{,}160 & \textbf{2784.3} & -- & -- & -- & -- & -- & -- & -- \\
 & 32 & 41{,}746 & \textbf{2801.5} & -- & -- & -- & -- & -- & -- & -- \\
\bottomrule
\end{tabular}
\end{table}

\paragraph{Overall comparison.}
Figure~\ref{fig:solver-comparison} summarises both benchmarks by the shifted geometric mean of the solve time (shift $10$\,s), normalised by that of the fastest solver in each panel. Timeouts, out-of-memory runs, runs whose relative objective error exceeds $10^{-3}$, and configurations
not attempted due to computational resource
limitations are assigned $3600$\,s in the aggregate comparison. On the PageRank instances, sPDLP with snapping solves all $18$ and is the fastest under this measure: single-threaded MOSEK is $1.69$ times slower, PDLP $5.44$ times and sPDLP without snapping $10.55$ times. Many runs without snapping either reach a resource limit or stop with walls short of their bounds. On the $68$ MIPLIB relaxations of Section~\ref{sec:experiments}, the three first-order variants are within about $13$ per cent of one another, and sPDLP with snapping is within $3$ per cent of PDLP, in agreement with the per-instance ratios of Section~\ref{sec:experiments} and Appendix~\ref{app:snapping_benchmark}; single-threaded interior-point MOSEK is about $13$ to $15$ times faster than all three. Snapping thus costs little on MIPLIB, and can help considerably in the controlled PageRank setting.

\paragraph{Scope.}
These experiments show that snapping can reduce the iteration count considerably when the walls appear in the PageRank constraints, \emph{with fixed steps and without preconditioning}. They do not show an advantage over PDLP in its default, preconditioned configuration. They concern the behaviour of the solver, not the effect of the transfers on the PageRank vector.

\section{The convergence score}
\label{app:score}

The dual of~\eqref{eq:LP} is $\max_{y,\lambda}\ q^\top y + l^\top\lambda^+ - u^\top\lambda^-$ subject to $c-K^\top y=\lambda$, $y_G\ge 0$ and $\lambda\in\Lambda$, where $y=(y_G,y_A)$ is split as in $Y$ and $\Lambda$ is the set of reduced costs compatible with the bounds. The three terms of the score~\eqref{eq:score} are, as in \citet{Applegate2021PracticalLL},
\[
r_{\mathrm{p}}^k := \frac{\bigl\|\bigl(Ax^k-b,\;(h-Gx^k)^+\bigr)\bigr\|_2}{1+\|q\|_2}, \qquad
r_{\mathrm{d}}^k := \frac{\|c - K^\top y^k - \lambda^k\|_2}{1+\|c\|_2},
\]
\[
g^k := \frac{|q^\top y^k + l^\top(\lambda^k)^+ - u^\top(\lambda^k)^- - c^\top x^k|}{1 + |q^\top y^k + l^\top(\lambda^k)^+ - u^\top(\lambda^k)^-| + |c^\top x^k|},
\qquad \lambda^k = \proj_\Lambda(c-K^\top y^k),
\]
the relative primal residual, the relative dual residual and the relative optimality gap.

\section{Proof of Theorem~\ref{thm:wall}}
\label{app:proof_wall}

Write $w_k$ for the $w$-coordinate of the $k$-th iterate and $W_k:=\max_{0\le i\le k}w_i$. Since $w$ enters no constraint, its component of $c-K^\top y$ equals $-\delta$ whatever $y$ is, so a PDHG step gives $w_{k+1}=\min\{1,\,w_k+\tau\delta\}\le W_k+\tau\delta$. A convex combination of earlier iterates has $w$-coordinate at most $W_k$. In either case $W_{k+1}\le W_k+\tau\delta$, and from $w_0=0$ induction gives $W_k\le k\tau\delta$. The target requires $w_k=1$, hence $k\tau\delta\ge1$ and $k\ge L$. The step-size condition gives $\tau<1/\|\hat K\|_2=1/\sqrt n$, so $L\ge 1/(\tau\delta)>\sqrt n/\delta$. This proves (i).

For (ii), each restart draws afresh from a fixed law, so the restart-free counts of successive epochs are independent and identically distributed, and Assumption~\ref{ass:N} holds. With probability one half the restart sets $w=1$; the PDHG step keeps $w=1$, and the block reaches its target within $f$ iterations, so $N\le f$ on that event. Hence $\phi(p)=\E[(1-p)^N]\ge\tfrac12(1-p)^f$, and Theorem~\ref{thm:renewal} gives
\[
T_p=\frac{(1-p)\bigl(1-\phi(p)\bigr)}{p\,\phi(p)}\le\frac{1-p}{p\,\phi(p)}\le\frac{2}{p\,(1-p)^{f-1}} .
\]
At $p=1/f$ this is $2f\,(1-1/f)^{-(f-1)}\le2\mathrm{e}f$, because $(1-1/f)^{f-1}\ge\mathrm{e}^{-1}$ for every $f\ge2$. \qed

\section{Connections to neural-network training}
\label{app:ml}

We discuss connections between
the multiple-wall model of Appendix~\ref{app:many-walls}
and gradient behaviour in neural networks,
as well as possible wall mechanisms in LP formulations
of ReLU-network training. These connections are
illustrative and do not extend Proposition~\ref{prop:manywalls}
to the coupled formulations.

\subsection{Vanishing and exploding gradients}

Each wall coordinate feels a gradient of only $\delta$, the analogue of a vanishing gradient \citep{pascanu2013difficulty}. The block's large norm $\sqrt n$ forces the step size below $1/\sqrt n$, much as large gradients in some directions cap the usable learning rate \citep{philipp2017exploding}. Neither slows the method alone: the plateau takes about $\|\hat K\|_2/\delta$ iterations, a ratio of the largest scale in the problem to the smallest gradient that still matters. In this convex setting that ratio is where the two gradient problems appear together; in a deep network both arise instead from products of Jacobians across layers, which a single linear programme has no reason to contain. The wall itself is a bound active at the optimum, as for a unit that should be inactive or a weight that should be zero. Scheme (C) is close to what deep-learning practice already does when it reinitialises units that have stopped learning: the dormant-neuron resets of \citet{sokar2023dormant}, continual backpropagation \citep{dohare2024loss}, and shrink-and-perturb \citep{ash2020warm}. Proposition~\ref{prop:manywalls} shows, in an idealised model, why such resets should act on individual units rather than on the whole network.

\subsection{Walls in linear programmes for training ReLU networks}

\citet{bienstock2023principled} construct, for a given architecture, loss and sample size, a single polytope with one face for each data set, such that optimising a linear function over the face of a data set trains the network to within an additive $\epsilon$. The polytope comes from a discretisation of the parameters into a binary problem of low treewidth. Its size is linear in the number of samples and exponential in the size of the architecture, which puts it beyond simplex and interior-point codes and within the reach of first-order methods such as PDLP. Two features of this formulation match the model. At an optimal vertex the lifted binary variables take the values $0$ and $1$, so that every such coordinate sits at a wall. Neighbouring points of the parameter grid differ in loss by an amount controlled by the grid spacing, so near-optimal alternatives should have small reduced costs, which are plateaus. A coordinate-wise restart then sets an undecided binary variable to a random bound, a randomised rounding of one variable. The resulting point may violate the coupling constraints, which PDHG tolerates, since its iterates need not be feasible.

A second source of walls is the convex reformulation of two-layer ReLU training by \citet{pilanci2020neural}, a group-sparse programme over the activation patterns of the data. With absolute loss and $\ell_1$ penalties in place of the squared loss and the group norm, it becomes a linear programme of the same shape; we use this analogue only as an illustration and do not claim that it is equivalent to training a network. At an optimum most pattern blocks are zero; a pattern that is nearly useful has a small reduced cost and drifts slowly to zero; and a coordinate-wise restart that sets a block to zero is the pruning of one neuron.

\subsection{What is not shown}

Proposition~\ref{prop:manywalls} is proved for decoupled walls, which presolve would remove, and under an assumption on the block. In the two neural-network LP formulations discussed above, the walls are coupled through the constraints, so the proposition does not apply to them directly; we have neither proved that their reduced costs are small nor measured the resulting iteration counts for these
formulations. The factor $\sqrt n$ comes from the operator norm and can be reduced by diagonal preconditioning, and a plateau delays termination under the score~\eqref{eq:score} only while $\delta$ is large compared with the tolerance.
The coupled-wall PageRank experiments
of Appendix~\ref{app:pagerank} test the snapping
mechanism in a different setting, but do not extend
Proposition~\ref{prop:manywalls} to these neural-network formulations.
What the proposition does establish is the direction a realistic test should take: restart schemes that act on individual coordinates, triggered by the sign test that identifies settled coordinates, and compared with global restarts on linear programmes of the kinds above.

\section{Proofs for Section~\ref{sec:exact}}
\label{app:proofs}

\subsection{Proof of Theorem~\ref{thm:renewal}}
\label{app:proof_renewal}

The formula is eq.~(20) of \citet{bonomo2021first}; we give a proof in our notation. Let $\Gamma$ be the index of the first restart, so that $\Gamma$ is geometric on $\{1,2,\dots\}$ with $\Prob[\Gamma>t]=(1-p)^t$, independent of $N$. Decompose the run into \emph{epochs}, an epoch being a maximal stretch of iterations between two consecutive restarts. By Assumption~\ref{ass:N} the epochs are independent and identically distributed; let $(N_i,\Gamma_i)_{i\ge1}$ be the corresponding independent copies.

Consider a single epoch. Under our accounting convention a restart is applied \emph{after} an iteration, so the epoch ends by absorption in $\cT$ if and only if no restart occurs after any of the iterations $1,\dots,N-1$, that is, if and only if $\Gamma\ge N$; and the number of iterations consumed by the epoch is $\ell:=\min(N,\Gamma)$ in either case. Hence the absorption probability is
\begin{equation}
  \phi_{\mathrm{abs}}:=\Prob[\Gamma\ge N]=\E[(1-p)^{N-1}]=\phi(p)/(1-p),
  \label{eq:phiabs}
\end{equation}
which is strictly positive because $\Prob[N<\infty]>0$ and at most $1$, with equality only in the degenerate case $\Prob[N=1]=1$, where $T_p\equiv1$. For the epoch length, summing the tail and using independence,
\[
  \E[\ell]
  =\sum_{t\ge0}\Prob[N>t]\,\Prob[\Gamma>t]
  =\sum_{t\ge0}\Prob[N>t]\,(1-p)^t
  =\frac{1-\E[(1-p)^{N}]}{p}
  =\frac{1-\phi(p)}{p},
\]
where the third equality is the identity $\sum_{t\ge0}u^t\,\Prob[N>t]=\bigl(1-\E[u^N]\bigr)/(1-u)$ for $u\in[0,1)$, obtained by exchanging the order of summation. In particular $\E[\ell]\le\E[\Gamma]=1/p<\infty$.

Let $J:=\min\{i\ge1: \Gamma_i\ge N_i\}$ be the index of the absorbing epoch. Then $J$ is geometric with success probability $\phi_{\mathrm{abs}}$, so $\E[J]=1/\phi_{\mathrm{abs}}<\infty$, and $\{J\ge i\}$ is determined by epochs $1,\dots,i-1$, hence independent of $\ell_i$. The total iteration count is $\sum_{i=1}^{J}\ell_i$, and Wald's identity gives
\[
  T_p=\E\Bigl[\sum_{i=1}^{J}\ell_i\Bigr]
     =\E[J]\,\E[\ell]
     =\frac{1-p}{\phi(p)}\cdot\frac{1-\phi(p)}{p},
\]
which is~\eqref{eq:renewal}, and is finite. Finally, $\phi(p)\to\Prob[N<\infty]$ and $(1-\phi(p))/p=\sum_{t\ge0}\Prob[N>t](1-p)^t\to\E[N]$ as $p\downarrow0$ by monotone convergence, which gives the last claim. \qed

\subsection{Two extreme cases}

\begin{remark}[Two extreme cases]\label{rem:extremes}
If $N\equiv t$ is deterministic then~\eqref{eq:renewal} telescopes: with $u=1-p$,
\begin{equation}
  T_p=\frac{u\,(1-u^{t})}{(1-u)\,u^{t}}
     =\sum_{k=0}^{t-1}u^{-k}\;=:\;\Tdet_p(t)\;\ge\;t ,
  \label{eq:detT}
\end{equation}
and $T_p$ is strictly increasing in $p$ for $t\ge2$: restarting a method whose running time is exactly predictable can only hurt, at every $p$. If instead $N$ is geometric, so that the restart-free method is memoryless, then~\eqref{eq:renewal} gives $T_p=\E[N]$ for all $p$, and~\eqref{eq:cvcriterion} holds with equality.
\end{remark}

\subsection{An arbitrary starting point}

Theorem~\ref{thm:renewal} assumes that the run begins at a fresh restart point; the version for a solve that begins at a user-supplied $z^0$ follows from the same one-epoch decomposition.

\begin{corollary}[Arbitrary starting point]\label{cor:start}
Let $N_0$ be the restart-free iteration count from $z^0 \notin\mathcal{T}$ and $\phi_0(p):=\E[(1-p)^{N_0}]$. Then the expected iteration count from $z^0$ is
\begin{equation}
  T_p(z^0)\;=\;\frac{1-\phi_0(p)}{p}
  \;+\;\Bigl(1-\frac{\phi_0(p)}{1-p}\Bigr)\,T_p .
  \label{eq:renewal_start}
\end{equation}
\end{corollary}
\begin{proof}
The first epoch consumes $\E[\min(N_0,\Gamma)]=(1-\phi_0(p))/p$ iterations and absorbs with probability $\phi_0(p)/(1-p)$, as in~\eqref{eq:phiabs}. Otherwise the process starts afresh from $\mu$ and contributes $T_p$.
\end{proof}

\subsection{Proof of Corollary~\ref{cor:cv}}
\label{app:proof_dispersion}

The criterion is~(25) of \citet{bonomo2021first}; we give a proof in our notation. Let $\psi(u):=\E[u^{N}]$ be the probability generating function of $N$, so that $\phi(p)=\psi(1-p)$. Since $\E[N^2]<\infty$, Abel's theorem gives $\psi(1)=1$, $\psi'(1^-)=\E[N]$ and $\psi''(1^-)=\E[N(N-1)]$, hence
\[
  \phi(p)=1-pm+p^{2}A+o(p^{2}),
  \qquad m:=\E[N],\quad
  A:=\tfrac12\E[N(N-1)]=\tfrac12\bigl(v+m^{2}-m\bigr),
\]
with $v:=\Var(N)$. Substituting into~\eqref{eq:renewal},
\[
  \frac{1-\phi(p)}{p}=m-pA+o(p),
  \qquad
  \frac{1}{\phi(p)}=1+pm+o(p),
\]
so that $\frac{1-\phi}{p\phi}=m+p\,(m^{2}-A)+o(p)$ and therefore
\[
  T_p=(1-p)\Bigl(m+p(m^{2}-A)\Bigr)+o(p)
     =m+p\bigl(m^{2}-m-A\bigr)+o(p)
     =m+\frac{p}{2}\bigl(m^{2}-m-v\bigr)+o(p).
\]
This is~\eqref{eq:cvderiv}, and both implications follow from the sign of the first-order term. \qed

\subsection{Proofs of Theorem~\ref{thm:mrl} and Corollary~\ref{cor:nbue}}
\label{app:proof_mrl}

Write $\bar F(j):=\Prob[N>j]$ for the tail of $N$, so that $\E[(N-j)^+]=\sum_{i\ge j}\bar F(i)$ and $m(j)=\E[(N-j)^+]/\bar F(j)$ whenever $\bar F(j)>0$. Throughout, a term $\bar F(j)\,m(j)$ with $\bar F(j)=0$ is read as $\E[(N-j)^+]=0$.

\paragraph{Proof of Theorem~\ref{thm:mrl}.}
Let $\Gamma$ be geometric on $\{1,2,\dots\}$ with $\Prob[\Gamma=j]=p(1-p)^{j-1}$, independent of $N$. From the proof of Theorem~\ref{thm:renewal}, $T_p=\E[\min(N,\Gamma)]/\Prob[\Gamma\ge N]$ with $\Prob[\Gamma\ge N]=\phi(p)/(1-p)$. Since $N=\min(N,\Gamma)+(N-\Gamma)^+$ and $m<\infty$,
\[
  \E[\min(N,\Gamma)]
  =m-\E[(N-\Gamma)^+]
  =m-\sum_{j\ge1}\Prob[\Gamma=j]\,\E[(N-j)^+]
  =m-\sum_{j\ge1}\Prob[\Gamma=j]\,\bar F(j)\,m(j).
\]
Also $\Prob[\Gamma<N]=\sum_{j\ge1}\Prob[\Gamma=j]\,\bar F(j)$, so $m=m\,\Prob[\Gamma\ge N]+\sum_{j\ge1}\Prob[\Gamma=j]\,\bar F(j)\,m$. Subtracting,
\[
  \E[\min(N,\Gamma)]-m\,\Prob[\Gamma\ge N]
  =\sum_{j\ge1}\Prob[\Gamma=j]\,\bar F(j)\,\bigl(m-m(j)\bigr),
\]
and dividing by $\Prob[\Gamma\ge N]=\phi(p)/(1-p)>0$ gives~\eqref{eq:mrl}. \qed

\paragraph{A memorylessness lemma.}
The strictness claim of Corollary~\ref{cor:nbue} rests on the following elementary fact.

\begin{lemma}\label{lem:memoryless}
Let $N$ take values in $\{1,2,\dots\}$ with $m=\E[N]<\infty$. Then $m(j)=m$ for every $j\ge0$ with $\bar F(j)>0$ if and only if $N$ is geometric with parameter $1/m$.
\end{lemma}
\begin{proof}
A geometric $N$ is memoryless, so $m(j)=m$. Conversely, suppose $m(j)=m$ whenever $\bar F(j)>0$, that is, $\sum_{i\ge j}\bar F(i)=m\,\bar F(j)$. Fix $j$ with $\bar F(j)>0$. If $\bar F(j+1)=0$ the identity reads $\bar F(j)=m\,\bar F(j)$, so $m=1$ and $N\equiv1$, which is geometric with parameter $1$. Otherwise subtract the identity at $j+1$ from the one at $j$ to obtain $\bar F(j)=m\bigl(\bar F(j)-\bar F(j+1)\bigr)$, that is, $\bar F(j+1)=(1-1/m)\,\bar F(j)$. Starting from $\bar F(0)=1$, induction gives $\bar F(j)=(1-1/m)^{j}$ for all $j$, the tail of a geometric law with parameter $1/m$.
\end{proof}

\paragraph{Proof of Corollary~\ref{cor:nbue}.}
(i) and (ii): every weight $\frac{1-p}{\phi(p)}\,p(1-p)^{j-1}\bar F(j)$ in~\eqref{eq:mrl} is non-negative, and strictly positive whenever $\bar F(j)>0$, so the sign of $T_p-T_0$ is that of the terms $m-m(j)$ when these all have the same sign. For strictness in (i): if $N$ is not geometric, Lemma~\ref{lem:memoryless} gives some $j\ge1$ with $\bar F(j)>0$ and $m(j)\ne m$; under NBUE this means $m(j)<m$, so that term is strictly positive and $T_p>T_0$ for every $p\in(0,1)$. If $N$ is geometric, Remark~\ref{rem:extremes} gives $T_p=T_0$.

(iii): For $j\ge0$ and $k\ge1$ with $\bar F(j)>0$,
\[
  \mathbb{P}[N>j+k\mid N>j]
=\prod_{t=j+1}^{j+k}(1-H(t)).
\]
If $H$ is non-decreasing, each factor $1-H(j+i)$, $i=1,\dots,k$, is non-increasing in $j$, so the conditional tail is non-increasing in $j$ for every $k$, and therefore $m(j)=\sum_{k\ge0}\Prob[N>j+k\mid N>j]$ is non-increasing in $j$. In particular $m(j)\le m(0)=m$, which is NBUE. This is the discrete form of the classical chain IFR $\Rightarrow$ DMRL $\Rightarrow$ NBUE \citep{barlow1975statistical}. \qed

\paragraph{Recovering Corollary~\ref{cor:cv} from Theorem~\ref{thm:mrl}.}
Assume $\E[N^2]<\infty$. Then $\sum_{j\ge1}\bar F(j)=m-1$ and $\sum_{j\ge1}\bar F(j)\,m(j)=\sum_{j\ge1}\E[(N-j)^+]=\E\bigl[\tfrac12N(N-1)\bigr]$, so $\sum_{j\ge1}\bar F(j)\,|m-m(j)|<\infty$. As $p\downarrow0$, $(1-p)/\phi(p)\to1$ and, by dominated convergence, $\sum_{j\ge1}(1-p)^{j-1}\bar F(j)\bigl(m-m(j)\bigr)\to\sum_{j\ge1}\bar F(j)\bigl(m-m(j)\bigr)=m(m-1)-\tfrac12\E[N(N-1)]$. Hence, with $v=\Var(N)$,
\[
  T_p-T_0=p\Bigl(m(m-1)-\tfrac12\bigl(v+m^2-m\bigr)\Bigr)+o(p)
         =\frac p2\bigl(m^2-m-v\bigr)+o(p),
\]
which is~\eqref{eq:cvderiv}.

\subsection{Proofs of Proposition~\ref{prop:bounded} and Corollary~\ref{cor:linear}}
\label{app:proof_bounded}

Let $\Gamma$ be as above. Since $t_-\le N\le t_+$ almost surely, $\min(t_-,\Gamma)\le\min(N,\Gamma)\le\min(t_+,\Gamma)$ pointwise and $\Prob[\Gamma\ge t_+]\le\Prob[\Gamma\ge N]\le\Prob[\Gamma\ge t_-]$, all three probabilities being positive. Hence
\[
  \frac{\E[\min(t_-,\Gamma)]}{\Prob[\Gamma\ge t_-]}
  \;\le\;
  \frac{\E[\min(N,\Gamma)]}{\Prob[\Gamma\ge N]}=T_p
  \;\le\;
  \frac{\E[\min(t_+,\Gamma)]}{\Prob[\Gamma\ge t_+]} .
\]
The outer expressions are the values of~\eqref{eq:renewal} for the deterministic counts $t_-$ and $t_+$, which Remark~\ref{rem:extremes} evaluates as $\Tdet_p(t_\pm)=\sum_{k=0}^{t_\pm-1}(1-p)^{-k}$. This is the sandwich in~\eqref{eq:bounded}. Since every summand of $\Tdet_p(t_-)$ is at least $1$, $T_p\ge t_-$ for every $p$, while $T_0=\E[N]\le t_+$; dividing gives the bound on the speed-up.

For the last claim, Popoviciu's inequality for a random variable supported in $[t_-,t_+]$ gives $\Var(N)\le(t_+-t_-)^2/4$, and $x\mapsto x^2-x$ is increasing on $[1,\infty)$ with $m\ge t_-\ge1$, so $m^2-m\ge t_-(t_--1)$. Under the stated condition, $\Var(N)\le(t_+-t_-)^2/4<t_-(t_--1)\le m^2-m$, that is, the inequality in~\eqref{eq:cvcriterion} is reversed and strict; $\E[N^2]<\infty$ holds trivially, so Corollary~\ref{cor:cv} gives $T_p>T_0$ for all sufficiently small $p>0$. \qed

\paragraph{Proof of Corollary~\ref{cor:linear}.}
The difference $t_+(\varepsilon)-t_-(\varepsilon)$ is bounded by $\Delta:=\log(\rho_+/\rho_-)/\log(1/\gamma)+t_0+1$ independently of $\varepsilon$, while $t_-(\varepsilon)\to\infty$ as $\varepsilon\downarrow0$. Claim (i) is the speed-up bound of Proposition~\ref{prop:bounded}, since $t_+/t_-\le1+\Delta/t_-\to1$. Claim (ii) is the last sentence of the proposition, whose condition $(t_+-t_-)^2<4t_-(t_--1)$ holds once $\Delta^2<4t_-(t_--1)$, that is, for all $\varepsilon$ below some $\varepsilon_0$. \qed

\subsection{Charging the restart}

\begin{remark}[Charging the restart]\label{rem:convention}
If each restart is charged one iteration, the same argument with $\Prob[\Gamma>N]$ in place of $\Prob[\Gamma\ge N]$ in~\eqref{eq:phiabs} yields $T_p=(1-\phi(p))/(p\,\phi(p))$ and the criterion $\CV(N)^2>1+1/\E[N]$. Restarting is then strictly harmful in the memoryless case, by exactly the wasted iteration. The structure of every statement in Section~\ref{sec:exact} is unaffected; the constants are not. As Remark~\ref{rem:preempt} notes, this convention yields the same formula as letting a restart pre-empt an absorption, which is the convention of Section~\ref{sec:threshold}.
\end{remark}

\section{Proofs for Section~\ref{sec:threshold}}
\label{app:proof_threshold}

\begin{remark}[A restart may pre-empt absorption]\label{rem:preempt}
The Bellman recursion~\eqref{eq:bellman} applies the restart action with probability $a$ at every state outside $\cT$, irrespective of where the accompanying iteration would have landed. Repeating the one-epoch argument of Appendix~\ref{app:proof_renewal} with $\Prob[\Gamma>N]=\phi(p)$ in place of $\Prob[\Gamma\ge N]$ gives, for a constant rate $p$,
\begin{equation}
  \tilde T_p=\frac{1-\phi(p)}{p\,\phi(p)}=\frac{T_p}{1-p},
  \qquad
  \tilde T_p(z^0)=\frac{1-\phi_0(p)}{p}+\bigl(1-\phi_0(p)\bigr)\tilde T_p
  \;\le\;\frac1p+\tilde T_p ,
  \label{eq:preempt}
\end{equation}
which is the expression of Remark~\ref{rem:convention}; the second formula is the counterpart of Corollary~\ref{cor:start} and gives the uniform bound~\eqref{eq:Tbar}. All statements of Section~\ref{sec:threshold} are qualitative and hold under either convention.
\end{remark}

\paragraph{Proof of Lemma~\ref{lem:bangbang}.}
Equation~\eqref{eq:bellman} is the Bellman equation of the dynamic programme: one iteration is always performed, after which, with probability $1-a$, the state is distributed as $P(\cdot\mid R)$ and, with probability $a$, as $\mu$. The map $a\mapsto Q^\star(R)+a\bigl(C^\star-Q^\star(R)\bigr)$ is affine, so its minimum over the compact set $\cA$ is attained at $\max\cA$ when the slope is negative, at $\min\cA$ when it is positive, and everywhere when it vanishes. Measurability of the optimal selector is elementary, since it is determined by the sign of $Q^\star(R)-C^\star$. \qed

\paragraph{Proof of Theorem~\ref{thm:threshold}.}
(i) Set $T^{(0)}\equiv0$ and, for $j\ge0$,
\[
  T^{(j+1)}(R):=\mathbb{1}\{R>\varepsilon\}
  \Bigl(1+\min_{a\in\cA}\bigl[(1-a)\,Q^{(j)}(R)+a\,C^{(j)}\bigr]\Bigr),
\]
where $Q^{(j)}(R)=\int T^{(j)}\,dP(\cdot\mid R)$ and $C^{(j)}=\int T^{(j)}\,d\mu$. Under Assumption~\ref{ass:ssp}, $T^{(j)}\uparrow T^\star$ pointwise, with $T^\star\le\bar T$ \citep[Ch.~8--9]{bertsekas1978stochastic}. We show by induction that each $T^{(j)}$ is non-negative and non-decreasing. This holds for $T^{(0)}$. If it holds for $T^{(j)}$, then $T^{(j)}$ is bounded and non-decreasing, so Assumption~\ref{ass:mono} makes $Q^{(j)}$ non-decreasing, while $C^{(j)}$ is a constant. By affineness the minimum over $\cA$ equals the minimum of the two values at $\underline a=\min\cA$ and $\bar a=\max\cA$, each a non-decreasing affine function of $Q^{(j)}(R)$ with slope $1-a\ge0$, and a minimum of non-decreasing functions is non-decreasing. Multiplying by $\mathbb{1}\{R>\varepsilon\}$ preserves monotonicity, because the indicator is non-decreasing and the other factor non-negative. Letting $j\to\infty$, $T^\star$ is non-decreasing, and by Assumption~\ref{ass:mono} so is $Q^\star$.

(ii) If $R\le R'$ and $Q^\star(R)>C^\star$, then $Q^\star(R')\ge Q^\star(R)>C^\star$ by (i); thus the switching set is an
up-set in $\mathcal S$, of the stated form
after intersection with $\mathcal S$.

(iii)
By Lemma~\ref{lem:bangbang}, the action
$\max\cA$ is optimal when $Q^\star(R)>C^\star$,
and $\min\cA$ is optimal otherwise. By part~(ii), this defines a two-level threshold policy, with the action at
$\beta^\star$, when $\beta^\star\in\cS$,
determined by the sign of
$Q^\star(\beta^\star)-C^\star$. A stationary policy attaining the minimum in~\eqref{eq:bellman} at every state is optimal for a dynamic programme with non-negative costs \citep[Ch.~8--9]{bertsekas1978stochastic}. \qed

\paragraph{Proof of Corollary~\ref{cor:plow}.}
The first claim follows from~\eqref{eq:betastar} and~\eqref{eq:optpolicy}: for $R<\beta^\star$, $Q^\star(R)\le C^\star$, so $\pi^\star(R)=\min\cA=0$. For the second, $T^\star$ is non-decreasing and vanishes on $\cT$, so $C^\star\ge \omega\,T^\star(\xi)\ge \omega$, where $T^\star(\xi)\ge1$ by~\eqref{eq:bellman} since $\xi>\varepsilon$. On the other side, $T^\star\le\bar T$ and $T^\star=0$ on $\cT$, so $Q^\star(R)\le\bar T\,\Prob[R_1>\varepsilon\mid R_0=R]$. By~\eqref{eq:plow_condition} this is $<\omega\,T^\star(\xi)\le C^\star$ on a right neighbourhood of $\varepsilon$, so $\beta^\star>\varepsilon$ by~\eqref{eq:betastar}. Since $T^\star(\xi)\ge1$, the condition $\limsup_{R\downarrow\varepsilon}\Prob[R_1>\varepsilon\mid R_0=R]<\omega/\bar T$ suffices. \qed

\begin{remark}[Interpretation and calibration of the threshold]\label{rem:calib}

By~\eqref{eq:betastar} and
Lemma~\ref{lem:bangbang},
$\beta^\star$ separates the states where $Q^\star(R)>C^\star$ from those where $Q^\star(R)\le C^\star$. Equality at $R=\beta^\star$, when
$\beta^\star\in\cS$, is not guaranteed without additional regularity. The cutoff $\kappa\varepsilon_{\mathrm{opt}}$ of~\eqref{eq:no_restart_region} plays the role of $\beta^\star$, although it is imposed on the score rather than on $R$, and no single global value can be right for every instance, since $\beta^\star$ depends on the instance through both $Q^\star$ and $C^\star$. Calibrating it exactly is circular, since $Q^\star$ and $C^\star$ are defined through $T^\star$. Estimates of $T_p$ for constant $p\in\cA\cap(0,1)$, obtained from pilot runs using Theorem~\ref{thm:renewal}, can guide constant-rate selection but do not by themselves determine $\beta^\star$.

\end{remark}

\begin{remark}[Why randomise at all?]\label{rem:randomise}

When $0,1\in\cA$,
Theorem~\ref{thm:threshold} admits an optimal
deterministic restart-or-continue policy:
restart when $Q^\star(R)>C^\star$ and continue
otherwise. Randomisation confers no benefit
in this idealised model, and the cap $p_{\max}<1$ and the factor $2^{0.6\lfloor\nu_k/d\rfloor}$ in~\eqref{eq:reff} bound the rate at which progress can be discarded. The situation changes once one recalls that the merit $M_k$ is a lossy statistic of the state $z^k$: a policy measurable with respect to $M_k$ is a restricted class in a partially observed problem, and there randomised memoryless policies can be strictly better than deterministic ones \citep{singh1994learning}. This offers one reason for the stochastic form of~\eqref{eq:restart_probability}; we have not shown that it applies here.
\end{remark}

\section{The diffusion surrogate}
\label{app:surrogate}

As a heuristic guide, and as a concrete model in which Assumption~\ref{ass:mono} can be verified, we consider the one-dimensional surrogate
\begin{equation}
dR_t = -V'(R_t)\,dt + \sqrt{2D}\,dW_t,
\qquad V(R) = \tfrac{1}{2}R^2,
\label{eq:SDE}
\end{equation}
where $D>0$ is an effective diffusion coefficient and $W$ a Brownian motion. The quadratic potential encodes mean reversion. With a state-dependent restart rate $\varrho(R)$ and restart law $\mu_R$, the mean first passage time of~\eqref{eq:SDE} satisfies the backward equation $D\,T''(R)-R\,T'(R)-\varrho(R)\bigl(T(R)-C(R)\bigr)=-1$ with $C(R)=\int T\,d\mu_R$ \citep{Evans_2011}, so that restarting at $R$ is beneficial precisely when $C(R)<T(R)$. We do not use the surrogate for quantitative predictions.

\begin{lemma}[The diffusion surrogate is monotone]\label{lem:diffmono}
Let $R_k$ evolve by the Euler--Maruyama discretisation of~\eqref{eq:SDE}, $R_{k+1}=R_k+\Delta t\,f(R_k)+\sqrt{2D\Delta t}\,w_k$, with $w_k$ independent and identically distributed and $f=-V'$ Lipschitz with constant $\mathrm{Lip}(f)$. If $\Delta t\,\mathrm{Lip}(f)\le 1$ then the transition kernel $P$ satisfies Assumption~\ref{ass:mono}. In particular this holds for $V(R)=\tfrac12R^2$ and any $\Delta t\le1$.
\end{lemma}
\begin{proof}
Couple two copies driven by the same noise, started at $R_0\le R_0'$. Then $R_1'-R_1=(R_0'-R_0)+\Delta t\bigl(f(R_0')-f(R_0)\bigr)\ge\bigl(1-\Delta t\,\mathrm{Lip}(f)\bigr)(R_0'-R_0)\ge0$. A monotone pathwise coupling is equivalent to stochastic monotonicity of $P$.
\end{proof}

Projecting the discretised state onto $[0,\infty)$
preserves stochastic monotonicity and leaves the
first passage time to the target unchanged for
$\varepsilon\geq0$.

The lemma concerns a kernel driven by additive noise independent of the state. It says nothing about the kernel induced on $R_k$ by one projected PDHG iteration~\eqref{eq:pdlp}, for which, as noted in Section~\ref{sec:threshold}, we know neither a proof nor a counterexample.

\begin{remark}[Condition~\eqref{eq:plow_condition} is needed]\label{rem:plow_needed}
In the surrogate at a fixed step size the one-step kernel does not concentrate into the target: for the kernel of Lemma~\ref{lem:diffmono}, $\Prob[R_1>\varepsilon\mid R_0=R]\to\Prob[\sqrt{2D\Delta t}\,w_0>\varepsilon\Delta t]$ as $R\downarrow\varepsilon$, a fixed positive constant. Whether $\beta^\star>\varepsilon$ then depends on the size of this constant relative to $\omega\,T^\star(\xi)/\bar T$, which is what~\eqref{eq:plow_condition} asks. The argument that restarts become useless near the target because $T(R)\to0$ as $R\downarrow\varepsilon$ is valid for the continuous-time surrogate and false in discrete time, where $T^\star\ge1$ outside the target.
\end{remark}

\section{Further computational results}
\label{app:table}

\paragraph{Hardware and tuning.}
All experiments were run on a 64-bit Windows machine with an 11th-generation Intel Core i5-1145G7 processor at 2.60\,GHz and 16\,GB of memory. The hyperparameters of Algorithm~\ref{alg:stochastic_restart} were chosen by a grid search on the fourteen instances \texttt{30n20b8}, \texttt{aflow40b}, \texttt{air04}, \texttt{danoint}, \texttt{dfn-gwin-UUM}, \texttt{eil33-2}, \texttt{eilB101}, \texttt{lectsched-4-obj}, \texttt{macrophage}, \texttt{mcsched}, \texttt{ran16x16}, \texttt{reblock67}, \texttt{rmatr100-p10} and \texttt{rmine6}, which are among the $68$ instances reported; the resulting single setting was then used for every instance.

\paragraph{Full results.}
Table~\ref{tab:main} lists every instance of the benchmark. On \texttt{map18}, sPDLP reaches the time limit on all three seeds while PDLP terminates in $9\,400$ iterations; on \texttt{rail02}, PDLP reaches the time limit at $50\,080$ iterations without terminating, while sPDLP terminates on one of the three seeds, in $75\,481$ iterations, and reaches the time limit on the other two.

Table~\ref{tab:summary} summarises it.

\begin{table}[H]
  \caption{Summary of Table~\ref{tab:main}: sPDLP against PDLP at $\varepsilon_{\mathrm{opt}}=10^{-4}$, three seeds. Ratios are sPDLP over PDLP, so values below one favour sPDLP.}
  \label{tab:summary}
  \centering
  \footnotesize
  \begin{tabular}{lrr}
    \toprule
    Among the 66 instances with a ratio & Iterations & Time \\
    \midrule
    \quad ratio $<0.5$ & 1 & 2 \\
    \quad ratio in $[0.5,1)$ & 26 & 31 \\
    \quad ratio in $[1,1.5)$ & 27 & 14 \\
    \quad ratio in $[1.5,2)$ & 5 & 7 \\
    \quad ratio $\ge2$ & 7 & 12 \\
    \midrule
    Geometric-mean ratio & 1.074 & 1.078 \\
    \bottomrule
  \end{tabular}
\end{table}

\begingroup
\small
\setlength{\tabcolsep}{5pt}
\begin{longtable}{lrrrr}
\caption{Full results on the Mittelmann benchmark at $\varepsilon_{\mathrm{opt}}=10^{-4}$: sPDLP (Algorithm~\ref{alg:stochastic_restart}, three seeds) against PDLP with its default adaptive normalised-duality-gap restart. Ratios are geometric means, over the seeds that reached the tolerance on both sides, of the sPDLP figure divided by the PDLP figure; values below one favour sPDLP. $n$ and $m$ are the numbers of variables and constraints after presolve.}
\label{tab:main}\\
\toprule
Instance & $n$ & $m$ & Iter.\ ratio & Time ratio \\
\midrule
\endfirsthead
\toprule
Instance & $n$ & $m$ & Iter.\ ratio & Time ratio \\
\midrule
\endhead
\midrule
\multicolumn{5}{r}{\textit{continued on next page}} \\
\endfoot
\bottomrule
\endlastfoot
\texttt{30n20b8} & 18,380 & 576 & 0.711 & 1.189 \\
\texttt{aflow40b} & 2,728 & 1,442 & 2.237 & 1.108 \\
\texttt{air04} & 8,904 & 823 & 1.311 & 3.943 \\
\texttt{app1-2} & 26,871 & 53,467 & 0.045 & 0.078 \\
\texttt{bab5} & 21,600 & 4,964 & 0.918 & 1.483 \\
\texttt{beasleyC3} & 2,500 & 1,750 & 0.847 & 1.780 \\
\texttt{bienst2} & 505 & 576 & 1.121 & 1.038 \\
\texttt{binkar10\_1} & 2,298 & 1,026 & 0.771 & 1.913 \\
\texttt{bley\_xl1} & 5,831 & 175,620 & 1.051 & 1.184 \\
\texttt{buildingenergy} & 154,978 & 277,594 & 0.668 & 0.754 \\
\texttt{csched010} & 1,758 & 351 & 1.111 & 1.846 \\
\texttt{danoint} & 521 & 664 & 0.966 & 4.022 \\
\texttt{datt256} & 262,144 & 11,077 & 1.221 & 2.205 \\
\texttt{dfn-gwin-UUM} & 938 & 158 & 0.612 & 3.314 \\
\texttt{eil33-2} & 4,516 & 32 & 1.622 & 2.516 \\
\texttt{eilB101} & 2,818 & 100 & 0.912 & 2.678 \\
\texttt{ex10} & 17,680 & 69,608 & 0.937 & 2.444 \\
\texttt{glass4} & 322 & 396 & 4.421 & 1.755 \\
\texttt{gmu-35-40} & 1,205 & 424 & 1.344 & 2.193 \\
\texttt{iis-100-0-cov} & 100 & 3,831 & 1.121 & 1.174 \\
\texttt{iis-bupa-cov} & 345 & 4,803 & 1.123 & 1.494 \\
\texttt{iis-pima-cov} & 768 & 7,201 & 1.206 & 1.672 \\
\texttt{irish-e} & 61,728 & 104,259 & 0.939 & 0.972 \\
\texttt{lectsched-4-obj} & 7,901 & 14,163 & 1.295 & 2.656 \\
\texttt{macrophage} & 2,260 & 3,164 & 1.528 & 2.361 \\
\texttt{map18} & 164,547 & 328,818 & time limit & time limit \\
\texttt{mcsched} & 1,747 & 2,107 & 2.575 & 2.388 \\
\texttt{mik-250-1-100-1} & 251 & 151 & 0.557 & 1.333 \\
\texttt{mine-166-5} & 830 & 8,429 & 2.246 & 1.121 \\
\texttt{mine-90-10} & 900 & 6,270 & 1.138 & 0.877 \\
\texttt{n3div36} & 22,120 & 4,484 & 1.102 & 0.902 \\
\texttt{neos-1109824} & 1,520 & 28,979 & 2.007 & 0.914 \\
\texttt{neos-1337307} & 2,840 & 5,687 & 0.703 & 0.557 \\
\texttt{neos-1396125} & 1,161 & 1,494 & 1.284 & 0.621 \\
\texttt{neos-1601936} & 4,446 & 3,131 & 0.925 & 0.730 \\
\texttt{neos-686190} & 3,660 & 3,664 & 2.466 & 0.672 \\
\texttt{neos-849702} & 1,737 & 1,041 & 1.253 & 0.575 \\
\texttt{neos-916792} & 1,474 & 1,909 & 2.712 & 1.439 \\
\texttt{neos-934278} & 23,123 & 11,495 & 0.818 & 0.532 \\
\texttt{neos13} & 1,827 & 20,852 & 0.844 & 0.607 \\
\texttt{neos18} & 3,312 & 11,402 & 0.738 & 0.701 \\
\texttt{net12} & 14,115 & 14,021 & 0.855 & 0.702 \\
\texttt{netdiversion} & 129,180 & 119,589 & 1.109 & 0.686 \\
\texttt{newdano} & 505 & 576 & 1.121 & 0.573 \\
\texttt{ns1208400} & 2,883 & 4,289 & 1.312 & 0.726 \\
\texttt{ns1688347} & 2,685 & 4,191 & 0.796 & 0.548 \\
\texttt{ns1830653} & 1,629 & 2,932 & 0.720 & 0.771 \\
\texttt{opm2-z7-s2} & 2,023 & 31,798 & 1.591 & 0.829 \\
\texttt{qap15} & 22,275 & 6,330 & 1.085 & 0.654 \\
\texttt{rail02} & 270,869 & 95,791 & time limit & time limit \\
\texttt{rail507} & 63,019 & 509 & 0.876 & 0.629 \\
\texttt{ran16x16} & 512 & 288 & 1.172 & 0.440 \\
\texttt{reblock67} & 670 & 2,523 & 0.931 & 1.080 \\
\texttt{rmatr100-p10} & 7,359 & 7,260 & 0.526 & 0.689 \\
\texttt{rmatr100-p5} & 8,784 & 8,685 & 0.847 & 0.678 \\
\texttt{rmine6} & 1,096 & 7,078 & 1.204 & 1.733 \\
\texttt{rocII-4-11} & 9,234 & 21,738 & 0.713 & 0.728 \\
\texttt{rococoC10-001000} & 3,117 & 1,293 & 1.012 & 0.767 \\
\texttt{satellites1-25} & 9,013 & 5,996 & 1.103 & 0.698 \\
\texttt{sp98ic} & 10,894 & 825 & 1.711 & 1.294 \\
\texttt{sp98ir} & 1,680 & 1,531 & 0.923 & 1.121 \\
\texttt{stp3d} & 204,880 & 159,488 & 1.690 & 1.489 \\
\texttt{support10} & 14,770 & 165,684 & 1.451 & 0.827 \\
\texttt{tanglegram2} & 4,714 & 8,980 & 1.120 & 0.782 \\
\texttt{timtab1} & 397 & 171 & 1.442 & 0.517 \\
\texttt{unitcal\_7} & 25,755 & 48,939 & 0.929 & 0.779 \\
\texttt{vpphard} & 51,471 & 47,280 & 1.375 & 1.941 \\
\texttt{zib54-UUE} & 5,150 & 1,809 & 1.046 & 3.152 \\
\midrule
Geometric mean over the 66 instances with a ratio & & & 1.074 & 1.078 \\
\end{longtable}

\endgroup

\subsection{Mittelmann results with snapping}
\label{app:snapping_benchmark}

We also evaluate sPDLP with snapping on the same
Mittelmann benchmark, using the same tolerance,
time limit, and random seeds.
Table~\ref{tab:full-results-snap} reports the full results.
Over the 67 instances with a ratio, the geometric-mean
ratios against PDLP are $1.133$ in iterations and
$1.006$ in wall-clock time. sPDLP with snapping
outperforms PDLP on 27 instances in iterations
and 36 instances in time.
\begingroup
\small
\setlength{\tabcolsep}{5pt}
\begin{longtable}{lrrrr}
\caption{Full results on the Mittelmann benchmark at $\varepsilon_{\mathrm{opt}} = 10^{-4}$: sPDLP with snapping (three seeds) against PDLP with its default adaptive normalised-duality-gap restart. Ratios are geometric means, over the seeds that reached the tolerance on both sides, of the sPDLP figure divided by the PDLP figure; values below one favour sPDLP. $n$ and $m$ are the numbers of variables and constraints after presolve.}
\label{tab:full-results-snap}\\
\toprule
Instance & $n$ & $m$ & Iter.\ ratio & Time ratio \\
\midrule
\endfirsthead
\toprule
Instance & $n$ & $m$ & Iter.\ ratio & Time ratio \\
\midrule
\endhead
\midrule
\multicolumn{5}{r}{\emph{continued on next page}} \\
\endfoot
\bottomrule
\endlastfoot
\texttt{30n20b8} & 18,380 & 576 & 0.975 & 0.820 \\
\texttt{aflow40b} & 2,728 & 1,442 & 1.540 & 1.009 \\
\texttt{air04} & 8,904 & 823 & 1.234 & 1.012 \\
\texttt{app1-2} & 26,871 & 53,467 & 0.176 & 0.153 \\
\texttt{bab5} & 21,600 & 4,964 & 0.795 & 0.778 \\
\texttt{beasleyC3} & 2,500 & 1,750 & 0.774 & 0.955 \\
\texttt{bienst2} & 505 & 576 & 1.242 & 0.868 \\
\texttt{binkar10\_1} & 2,298 & 1,026 & 0.767 & 0.807 \\
\texttt{bley\_xl1} & 5,831 & 175,620 & 1.602 & 1.134 \\
\texttt{buildingenergy} & 154,978 & 277,594 & 0.874 & 0.775 \\
\texttt{csched010} & 1,758 & 351 & 1.345 & 1.108 \\
\texttt{danoint} & 521 & 664 & 0.891 & 0.910 \\
\texttt{datt256} & 262,144 & 11,077 & 1.229 & 1.043 \\
\texttt{dfn-gwin-UUM} & 938 & 158 & 0.545 & 0.916 \\
\texttt{eil33-2} & 4,516 & 32 & 1.346 & 2.103 \\
\texttt{eilB101} & 2,818 & 100 & 0.918 & 0.923 \\
\texttt{ex10} & 17,680 & 69,608 & 0.945 & 0.911 \\
\texttt{glass4} & 322 & 396 & 4.141 & 1.012 \\
\texttt{gmu-35-40} & 1,205 & 424 & 1.458 & 1.064 \\
\texttt{iis-100-0-cov} & 100 & 3,831 & 1.260 & 0.990 \\
\texttt{iis-bupa-cov} & 345 & 4,803 & 1.116 & 1.078 \\
\texttt{iis-pima-cov} & 768 & 7,201 & 1.175 & 1.093 \\
\texttt{irish-e} & 61,728 & 104,259 & 0.997 & 0.740 \\
\texttt{lectsched-4-obj} & 7,901 & 14,163 & 1.165 & 0.958 \\
\texttt{macrophage} & 2,260 & 3,164 & 1.387 & 0.996 \\
\texttt{map18} & 164,547 & 328,818 & 4.109 & 2.235 \\
\texttt{mcsched} & 1,747 & 2,107 & 2.526 & 1.769 \\
\texttt{mik-250-1-100-1} & 251 & 151 & 0.569 & 1.219 \\
\texttt{mine-166-5} & 830 & 8,429 & 2.103 & 0.987 \\
\texttt{mine-90-10} & 900 & 6,270 & 0.936 & 1.111 \\
\texttt{n3div36} & 22,120 & 4,484 & 1.256 & 1.203 \\
\texttt{neos-1109824} & 1,520 & 28,979 & 2.065 & 0.982 \\
\texttt{neos-1337307} & 2,840 & 5,687 & 0.682 & 0.826 \\
\texttt{neos-1396125} & 1,161 & 1,494 & 1.295 & 0.800 \\
\texttt{neos-1601936} & 4,446 & 3,131 & 0.855 & 0.911 \\
\texttt{neos-686190} & 3,660 & 3,664 & 1.938 & 1.193 \\
\texttt{neos-849702} & 1,737 & 1,041 & 1.198 & 0.998 \\
\texttt{neos-916792} & 1,474 & 1,909 & 3.242 & 2.254 \\
\texttt{neos-934278} & 23,123 & 11,495 & 0.729 & 1.360 \\
\texttt{neos13} & 1,827 & 20,852 & 0.785 & 0.768 \\
\texttt{neos18} & 3,312 & 11,402 & 0.789 & 0.971 \\
\texttt{net12} & 14,115 & 14,021 & 1.027 & 0.869 \\
\texttt{netdiversion} & 129,180 & 119,589 & 1.094 & 0.815 \\
\texttt{newdano} & 505 & 576 & 1.242 & 0.945 \\
\texttt{ns1208400} & 2,883 & 4,289 & 1.351 & 1.030 \\
\texttt{ns1688347} & 2,685 & 4,191 & 0.854 & 0.652 \\
\texttt{ns1830653} & 1,629 & 2,932 & 0.761 & 1.131 \\
\texttt{opm2-z7-s2} & 2,023 & 31,798 & 1.810 & 1.988 \\
\texttt{qap15} & 22,275 & 6,330 & 1.081 & 0.815 \\
\texttt{rail02} & 270,869 & 95,791 & time limit & time limit \\
\texttt{rail507} & 63,019 & 509 & 0.807 & 0.997 \\
\texttt{ran16x16} & 512 & 288 & 1.004 & 0.696 \\
\texttt{reblock67} & 670 & 2,523 & 0.903 & 1.193 \\
\texttt{rmatr100-p10} & 7,359 & 7,260 & 0.556 & 0.817 \\
\texttt{rmatr100-p5} & 8,784 & 8,685 & 0.753 & 1.026 \\
\texttt{rmine6} & 1,096 & 7,078 & 1.217 & 1.201 \\
\texttt{rocII-4-11} & 9,234 & 21,738 & 0.662 & 0.687 \\
\texttt{rococoC10-001000} & 3,117 & 1,293 & 0.916 & 0.880 \\
\texttt{satellites1-25} & 9,013 & 5,996 & 2.203 & 1.355 \\
\texttt{sp98ic} & 10,894 & 825 & 1.792 & 1.304 \\
\texttt{sp98ir} & 1,680 & 1,531 & 1.600 & 0.992 \\
\texttt{stp3d} & 204,880 & 159,488 & 1.443 & 1.003 \\
\texttt{support10} & 14,770 & 165,684 & 1.629 & 1.373 \\
\texttt{tanglegram2} & 4,714 & 8,980 & 1.179 & 1.366 \\
\texttt{timtab1} & 397 & 171 & 1.172 & 0.893 \\
\texttt{unitcal\_7} & 25,755 & 48,939 & 1.041 & 0.887 \\
\texttt{vpphard} & 51,471 & 47,280 & 1.260 & 1.480 \\
\texttt{zib54-UUE} & 5,150 & 1,809 & 0.955 & 1.033 \\
\midrule
\multicolumn{3}{l}{Geometric mean over the 67 instances with a ratio} & 1.133 & 1.006 \\
Geometric mean over the 66 common instances
& & & 1.111 & 0.994 \\
\end{longtable}
\endgroup
\paragraph{The instances map18 and rail02.}
On \texttt{map18}, sPDLP with snapping terminates on two
of the three seeds, whereas sPDLP without snapping reaches
the time limit on all three seeds. However, PDLP remains
faster on this instance. On \texttt{rail02}, PDLP reaches
the time limit, whereas sPDLP with snapping terminates
on two of the three seeds, compared with one seed without
snapping. The reported ratios are conditional on successful
runs and do not account for unsuccessful seeds.
For direct comparison with Table~\ref{tab:main}, Table~\ref{tab:full-results-snap} also reports
geometric-mean ratios over the common 66 instances.

\subsection{Controlled experiments on \texttt{danoint}}
\label{app:danoint}

\subsubsection{Renewal-formula test}
\label{app:danoint-renewal}

In this experiment, restart decisions
are made every 40 PDLP iterations, at the termination
checks. We therefore measure completion times in checks
rather than iterations, and $p$ is the restart probability
per check when applying the renewal formula~\eqref{eq:renewal}.

The fifteen restart-free counts behind the test of Section~\ref{sec:experiments} are
\[
N = 172,\,170,\,172,\,174,\,170,\,171,\,176,\,170,\,169,\,169,\,175,\,170,\,172,\,173,\,171 .
\]
Table~\ref{tab:danoint} gives the plug-in curve of~\eqref{eq:plugin} computed from them, together with the mean hitting time of independent runs of restarted PDHG at four values of $p$. The renewal formula predicts a monotone increase over the whole range, with its minimum at $p=0$. 

The observed means are consistent with a slowdown
at larger restart rates: at $p=0.0005$ the difference
from the control is not significant, at $p=0.001$
it is borderline, and at $p=0.01$ the slowdown
is clear but smaller than predicted. A likely reason for the gap is that coordinates without two finite bounds keep their pre-restart values, so post-restart states are not independent draws from $\mu$ and retain part of the progress made; in a replay experiment the actual post-restart states reached the target in $154.8$ checks on average against $171.6$ for fresh draws from $\mu$. Here $\mu$ denotes the calibration law used to generate restart-free initial
states: each primal coordinate with two finite bounds is sampled independently
and uniformly over its interval, coordinates without two finite bounds are set
to $0$, and the dual variable is initialised at $0$ (followed by projection).
The observed means at $p\in\{0,5\times10^{-4},10^{-3}\}$ are based on $15$
independent runs each, while the result at $p=10^{-2}$ uses $50$ independent
runs.

\begin{table}[H]
\caption{\texttt{danoint}: the plug-in curve $\hat T_p$ of~\eqref{eq:plugin} from fifteen restart-free counts, against the mean hitting time of independent runs of PDHG restarted at constant probability $p$, given as mean $\pm$ standard error with a $95$ per cent confidence interval.}
\label{tab:danoint}
\centering
\small
\begin{tabular}{ccccc}
\toprule
$p$ & $\hat\phi(p)$ & Predicted $\hat T_p$ & Observed mean $\pm$ s.e. & $95$\% interval \\
\midrule
0      & 1.0000 & 171.6  & $174.20\pm1.28$  & $[171.69,176.71]$ \\
0.0005 & 0.9178 & 179.1  & $172.27\pm0.83$  & $[170.65,173.88]$ \\
0.001  & 0.8422 & 187.1  & $198.40\pm12.30$ & $[174.29,222.51]$ \\
0.002  & 0.7093 & 204.6  & -- & -- \\
0.005  & 0.4231 & 271.3  & -- & -- \\
0.01   & 0.1783 & 456.3  & $361.20\pm28.07$ & $[306.18,416.22]$ \\
0.02   & 0.0312 & 1519.3 & -- & -- \\
\bottomrule
\end{tabular}
\end{table}

\subsubsection{Empirical NBUE test}
\label{app:danoint-nbue}

To examine the dichotomy of
Corollary~\ref{cor:nbue}, we extended the original 15 restart-free
counts to 200 independent draws from the same
calibration law $\mu$. The first 15 counts were
reproduced exactly, and no run was censored.
The 200 counts range from 165 to 184 checks,
with mean $\bar N=172.8$, standard deviation
$3.4$, and squared coefficient of variation
$\widehat{\mathrm{CV}}(N)^2\approx3.9\cdot10^{-4}$.

\begin{figure}[H]
\centering
\includegraphics[width=0.65\linewidth]{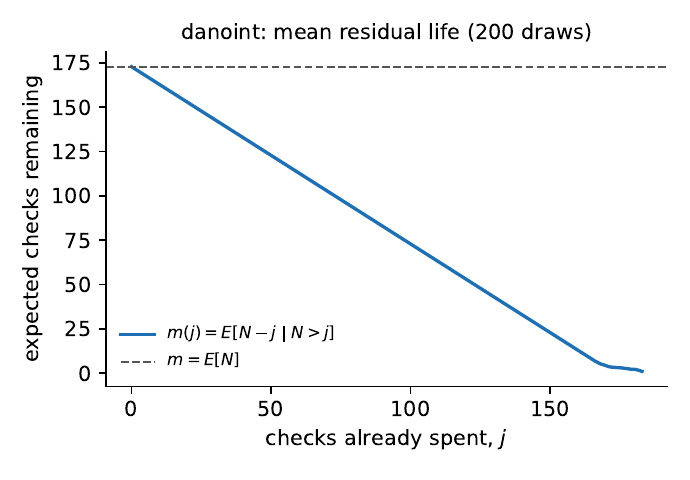}
\caption{Empirical mean residual
life $\widehat m(j)$ on \texttt{danoint}, computed
from 200 independent restart-free counts.
The dashed line indicates the restart-free
mean $\overline N=172.8$ checks.}
\label{fig:danoint-mrl}
\end{figure}

The empirical mean residual life
$\widehat m(j)$ satisfies $\widehat m(j)<\bar N$
for every $j\ge1$ with surviving observations.
Indeed, for $1\le j<165$, all sampled counts
exceed $j$, giving $\widehat m(j)=172.8-j$.
For $165\le j<184$, the remaining count is
at most $184-j\le19<172.8$. Thus, the
empirical distribution is NBUE, corresponding
to case~(i) of Corollary~\ref{cor:nbue}.
Figure~\ref{fig:danoint-mrl}
illustrates this result by comparing the
empirical mean residual life $\widehat{m}(j)$
with the restart-free sample mean $\bar{N}$.

The plug-in renewal predictions
from the 200 counts are 172.8, 180.5, 188.6,
and 463.1 checks for $p=0$, $0.0005$, $0.001$,
and $0.01$, respectively. They differ by less
than 3\% from the predictions based on the original 15 counts in Table~\ref{tab:danoint}.
Under the independent-restart
model, the empirical NBUE property predicts no
benefit from geometric restarting, consistent
with the observed slowdown at $p=0.01$
(from $174.2$ to $361.2$ checks). The NBUE conclusion applies
to the empirical distribution and does not
establish the condition for the unknown
restart-free count distribution.

We have not observed an NWUE
example and do not claim experimental
verification of Corollary~\ref{cor:nbue}(ii).

\subsubsection{Stochastic monotonicity diagnostic}
\label{app:danoint-monotonicity}

To investigate stochastic monotonicity,
we ran 50 independent restart-free PDLP trajectories
on \texttt{danoint}, initialised from independent draws
from $\mu$. We recorded $R_k$ at every iteration in
the solver's scaled space, using the iteration-dependent
step sizes $\tau_k=\eta_k/\omega_k$ and
$\sigma_k=\eta_k\omega_k$. We considered two horizons:
$h=1$ (one iteration) and $h=40$ (one termination check).

We partitioned the pooled values of $R_k$
into 20 equal-frequency bins. For each pair $b<b'$,
we estimated the conditional distribution functions $F_b(t)=\mathbb{P}(R_{k+h}\le t
\mid R_k\text{ belongs to bin }b)$ and $F_{b'}(t)$,
defined analogously for bin $b'$.
We then measured the largest violation of stochastic
ordering,
\begin{equation}
V_{b,b'}=\sup_t\bigl(F_{b'}(t)-F_b(t)\bigr)_+
\label{eq:CDFs_Violation}
\end{equation}
Here, $t$ is a threshold for
the future progress value $R_{k+h}$;
the supremum is evaluated over its observed values,
rather than at a preselected threshold.

For each horizon $h$, we formed
pairs $(R_k,R_{k+h})$ from the 50 trajectories;
these overlapping pairs are not independent.
For each bin pair, we use the descriptive KS
reference threshold
$1.36\sqrt{(n_b+n_{b'})/(n_bn_{b'})}$,
where $n_b$ and $n_{b'}$ are the respective
bin sample sizes.

\begin{table}[H]
\centering
\caption{Empirical stochastic-ordering
diagnostic on \texttt{danoint}: 50 restart-free trajectories
and 20 quantile bins. The KS reference is descriptive,
not a formal significance test, because observations
within trajectories are dependent.}
\label{tab:danoint-mono}
\begin{tabular}{lrrcc}
\toprule
& & Largest &
\multicolumn{2}{c}{Bin pairs above KS reference}\\
\cmidrule(lr){4-5}
Horizon & Observed pairs & violation & Adjacent & All\\
\midrule
$h=1$  & 324,785 & 0.042 & 2/19  & 3/190\\
$h=40$ & 322,835 & 0.424 & 11/19 & 118/190\\
\bottomrule
\end{tabular}
\end{table}

\begin{figure}[htbp]
\centering
\includegraphics[width=\linewidth]{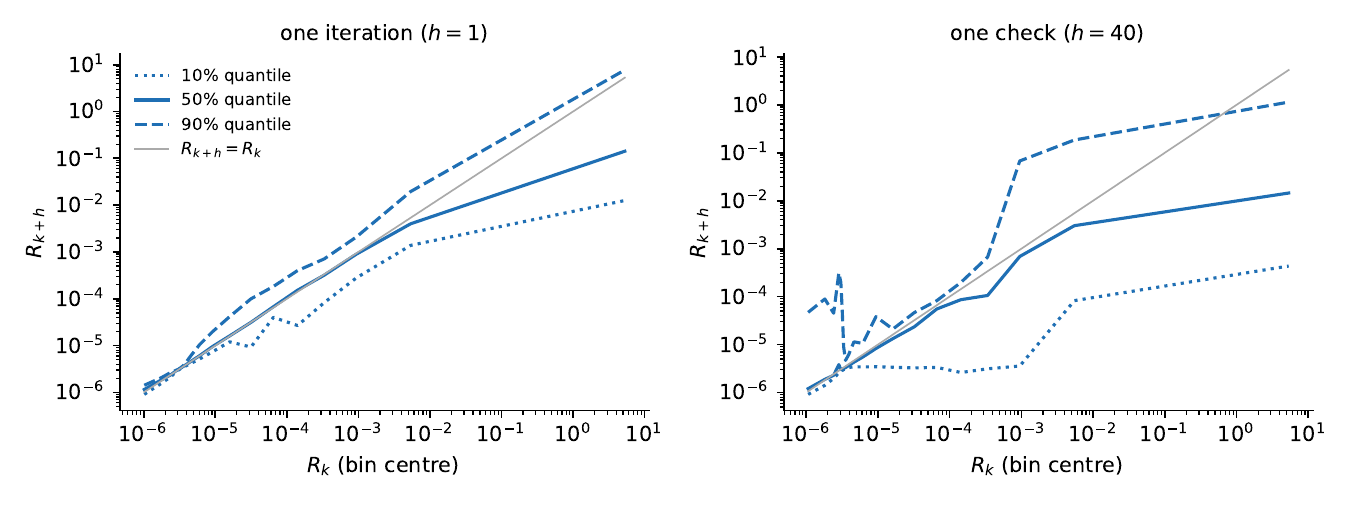}
\caption{Conditional 10\%, 50\%,
and 90\% quantiles of $R_{k+h}$ in each quantile bin
of $R_k$ on \texttt{danoint}, on logarithmic axes.
The grey line denotes $R_{k+h}=R_k$.
Left: $h=1$. Right: $h=40$.
Stochastic monotonicity predicts non-decreasing
conditional quantile curves.}
\label{fig:danoint-mono}
\end{figure}

Figure~\ref{fig:danoint-mono}
complements Table~\ref{tab:danoint-mono} by showing the conditional
10\%, 50\%, and 90\% quantiles of $R_{k+h}$
across bins of $R_k$ for both horizons.

At $h=1$, the observed departures
from stochastic ordering are small:
the largest violation is 0.042, and only 3 of 190 bin
pairs exceed the KS reference. At $h=40$ (40 iterations,
the interval between termination checks), the violations
are substantially larger, with a maximum of 0.424
and 118 of 190 bin pairs above the reference.
The conditional quantiles also show pronounced
non-monotonic behaviour at small $R_k$.

To check whether this result is caused
by excluding windows crossing the termination target,
we repeated the diagnostic with the target treated
as absorbing ($R=0$ after the first crossing).
At $h=40$, the maximum violation remains 0.419,
with 119 of 190 bin pairs above the reference.

If the reduced process were exactly
time-homogeneous and Markovian with a stochastically
monotone one-step kernel, its $h$-step kernel would
also be monotone. The difference between the two horizons therefore suggests limitations of the
reduced-state approximation. The experiment does not distinguish violations of exact
one-step monotonicity from non-Markovian or
time-varying effects, nor rule out estimation effects.
Assumption~\ref{ass:mono} remains unverified for projected PDHG,
and Theorem~\ref{thm:threshold} does not provide an optimality
guarantee for the practical sPDLP rule.

\subsection{Per-instance ratios and convergence profiles}
\label{app:profiles}

\begin{figure}[h]
  \centering
  \includegraphics[width=\textwidth]{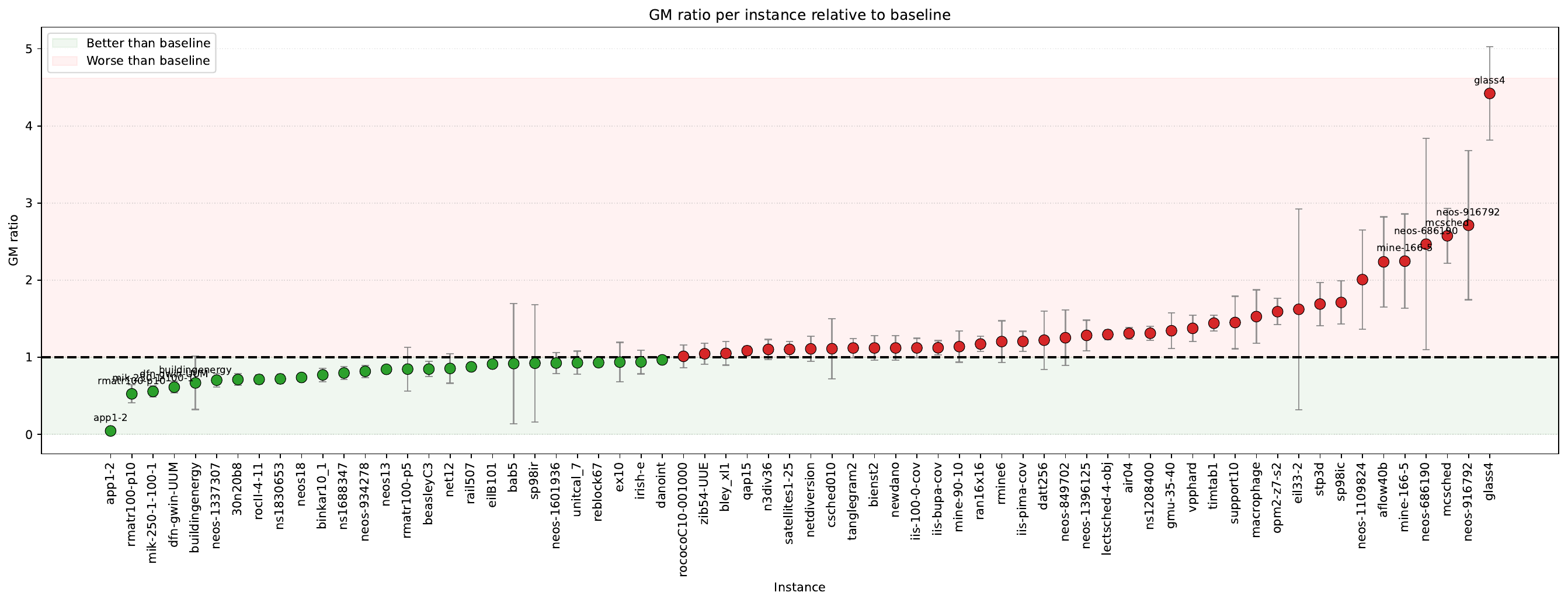}
  \caption{Per-instance geometric-mean iteration ratio of sPDLP to PDLP at $\varepsilon_{\mathrm{opt}}=10^{-4}$, instances ordered by ratio. The dashed line marks parity; points below it are instances on which sPDLP needs fewer iterations. Error bars show the range across the three seeds.}
  \label{fig:gmratio}
\end{figure}

Figure~\ref{fig:gmratio} shows the per-instance iteration ratios of Table~\ref{tab:main} in increasing order, with the range across the three seeds. The variability across seeds is small on most instances and substantial on a few.

\begin{figure}[h]
  \centering
  \includegraphics[width=\textwidth]{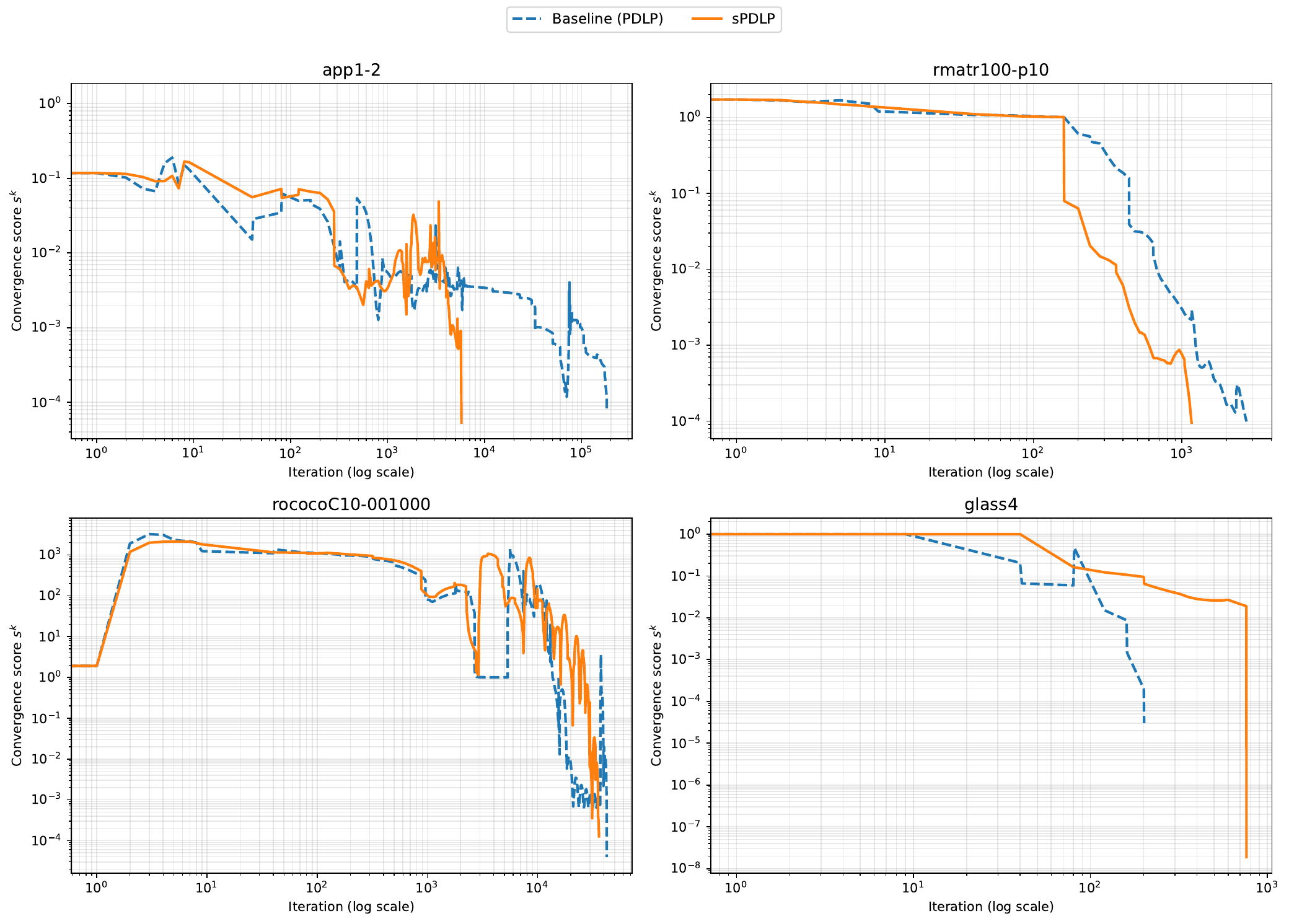}
  \caption{Convergence score $s^k$ against iteration count on four instances, both axes logarithmic. Dashed blue curves are PDLP and solid orange curves sPDLP.}
  \label{fig:convergence}
\end{figure}

Figure~\ref{fig:convergence} plots the convergence score $s^k$ against the iteration count for four instances. On \texttt{app1-2} and \texttt{rmatr100-p10}, sPDLP reaches small scores much earlier than PDLP, in keeping with their iteration ratios of $0.045$ and $0.526$. On \texttt{rococoC10-001000} the two methods behave alike, with a ratio of $1.012$. On \texttt{glass4}, PDLP reaches the target far earlier, while sPDLP stays at a larger score for many more iterations, in keeping with its ratio of $4.421$.

\subsection{Performance by instance size and sparsity}
\label{app:structure}

We ask whether the size or the sparsity of an instance helps to explain the relative performance of sPDLP. For each of the $66$ instances with a ratio we take its iteration ratio, the geometric mean over seeds of the sPDLP count divided by the PDLP count. Instances are grouped into terciles by the number of variables and by the density of the constraint matrix, the number of its nonzeros divided by the number of its entries; ties in density produce unequal tercile sizes. Table~\ref{tab:structural-performance} shows a tendency towards more frequent improvements on larger and sparser instances: the large-and-sparse group has a win rate of $56.3$ per cent against $33.3$ per cent for the small-and-dense group. The geometric mean of the large-and-sparse group is sensitive to \texttt{app1-2}, whose ratio is $0.045$; removing it moves the group's geometric mean from $0.836$ to $1.016$, while the win rate falls only from $56.3$ to $53.3$ per cent. Neither size nor sparsity yields a reliable rule for selecting instances, and the evidence does not establish an advantage in this regime.

\begin{table}[h]
\centering
\caption{Iteration performance by structural subgroup. The geometric-mean ratio is taken over the iteration ratios of the group; values below one favour sPDLP.}
\label{tab:structural-performance}
\small
\begin{tabular}{lrrr}
\toprule
Subgroup & Instances & Win rate & Geometric-mean ratio \\
\midrule
Small                    & 22 & 27.3\% & 1.236 \\
Large                    & 22 & 54.5\% & 0.891 \\
Dense                    & 23 & 34.8\% & 1.170 \\
Sparse                   & 22 & 54.5\% & 0.850 \\
\midrule
Small and dense          & 15 & 33.3\% & 1.181 \\
Large and sparse         & 16 & 56.3\% & 0.836 \\
\quad excluding \texttt{app1-2} & 15 & 53.3\% & 1.016 \\
\bottomrule
\end{tabular}
\end{table}

\end{document}